\documentclass{amsart}
\usepackage{amsfonts}

\newtheorem{theorem}{Theorem}
\theoremstyle{plain}

\newtheorem{conjecture}{Conjecture}
\newtheorem{corollary}{Corollary}

\newtheorem{lemma}{Lemma}

\newtheorem{remark}{Remark}

\numberwithin{equation}{section}

\begin{document}
\title[Selmer Groups]{Selmer Groups in Twist Families of Elliptic Curves}
\author{Ilker Inam}
\address[]{Uludag University, Faculty of Art and Science, Department of Mathematics\\
Gorukle, Bursa-Turkey}

\email[]{inam@uludag.edu.tr, ilker.inam@gmail.com}

\subjclass[2000]{14H52, 11G40, 14G10} \keywords{Elliptic curves,
Birch, Swinnerton-Dyer Conjecture, Zeta-functions and related
questions}

\begin{abstract}
The aim of this article is to give some numerical data related to
the order of the Selmer groups in twist families of elliptic
curves. To do this we assume the Birch and Swinnerton-Dyer
conjecture is true and we use a
celebrated theorem of Waldspurger to get a fast algorithm to compute $%
L_{E}(1)$. Having an extensive amount of data we compare the
distribution of the order of the Selmer groups by functions of
type $\alpha \frac{(\log \log (X))^{1+\varepsilon }}{\log (X)}$
with $\varepsilon $ small. We discuss how the "best choice" of
$\alpha $ is depending on the conductor of the chosen elliptic
curves and the congruence classes of twist factors.
\end{abstract}

\maketitle

\section{Introduction}

The purpose of this article is to give some numerical data related
to the order of the Selmer groups in twist families of elliptic
curves. This article is of experimental type. It would be most
interesting to give theoretical explanations for the statistical
observations we make.

Till today there is no algorithm that computes the order of the
Selmer group of a random elliptic curve defined over $\mathbb{Q}$
and even assuming the Conjecture of Birch and Swinnerton-Dyer it
is a hard computational problem to determine this order.

The situation will become easier if we restrict ourselves to twist
families of some specific elliptic curves. Working over
$\mathbb{Q}$ we can use the theory of modular forms to get an
analytic function for the $L$-series of the discussed curves.

\emph{Assuming }the Birch and Swinnerton-Dyer conjecture is true,
we are able to exploit a result of Waldspurger, which has a
crucial role in this article. It yields an efficient way to
compute the order of the Selmer groups in twist families of
elliptic curves, if one can find an eigenform of weight $3/2$
attached to the elliptic curve via the Shimura-Shintani lift.
Examples for this together with an explanation of how to apply
Waldspurger's result are discussed in \cite{ABF} (see Section
$4$). We use these examples and compute the orders of the Selmer
groups of twists of these curves up to $D\leq 10^{7}$. To do this
one begins with a curve $E$ and compares the order of the Selmer
groups of two twisted elliptic curves $E_{D_{0}}$ and $E_{D_{1}}$
with twist factors $D_{0}$ and $D_{1}$ in the same quadratic
congruence class modulo $4.N_{E}$ where $N_{E}$ is the conductor
of $E$. If one chooses the $E_{D_{0}}$ with $D_{0}$ small then its
Selmer group can be computed rather easily. So one can compute the
order of the Selmer groups for the elliptic curve $E_{D_{1}}$ by a
fast computation described in Subsection $3.3.$ After these
computations and with many data, it is a natural question to study
the distribution of members in twist families for which the Selmer
groups have the same order, say $k$ times the order of the torsion
of $E$ and to find simple functions that approximate this
distribution. Through the article, we only interested in $k$ in
order to compare the Selmer groups of different elliptic curves.
In this article we give numerical evidence that only constants
have to be changed for different twist families. We are interested
in twisted elliptic curves which have rank zero, but one has to be
careful about the cases where twisted elliptic curves have
(analytic) positive rank. We define $k=0$ to mean that the
corresponding twisted elliptic curve has positive (analytic) rank.
In this case the torsion subgroup doesn't play any role by
definition.

\subsection{Overview}

In Section $2$, we present some necessary definitions. The
notation used in the article is introduced. Section $3$ consists
of four subsections. In the first subsection, a statement of
Waldspurger's Theorem which plays a pivotal role in the article is
given. In Subsection $3.2.$, we describe how to compute
$d(n,n_{0})$. Proof is given which can be deduced from some
well-known facts. In Subsection $3.3.$ we describe the algorithm
to compute the order of the Selmer groups in twist families of
elliptic curves. Furthermore, the approximation function is
introduced in this subsection. We take the quotients of the
distribution functions and formulate a conjecture. Finally in
Section $4$, we give examples of our numerical results and in
particular tables listing constants $\alpha $ occurring in the
approximating functions.

Lastly we plot a graph showing the behavior of the distribution
function and the approximating function.

\subsection{Acknowledgements}

This article was partly written during my visit at the Institut
f\"{u}r Experimentelle Mathematik in Universit\"{a}t
Duisburg-Essen. I wish to express my gratitude for the support and
warm hospitality by this institution which made the visit a very
pleasant one and especially Prof.Dr.Gerhard Frey who suggested
this nice problem and made valuable comments and important
improvements on this article. Also I would like to thank
Prof.Dr.Gabor Wiese who made comments on an early version. This
article has grown out of my PhD thesis. This
article is supported by the The Scientific and Technological
Research Council of Turkey (TUBITAK) Research Project, Project No:
107T311. I wish to thank the referees for their helpful
suggestions.

\section{Background Material}

Let $E/\mathbb{Q}$ be an elliptic curve and assume that $D$ is a
square-free integer. With $E_{D}$ we denote the \emph{quadratic
twist }of $E$ with $D$. For $E$ given in "short" Weierstrass form
\begin{equation*}
y^{2}=x^{3}-g_{2}x-g_{3}.
\end{equation*}
$E_{D}$ is given by
\begin{equation*}
y^{2}=x^{3}-g_{2}D^{2}x-g_{3}D^{3}.
\end{equation*}
$E_{D}$ is the elliptic curve defined over $\mathbb{Q}$ isomorphic
to $E$ over $\mathbb{Q}(\sqrt{D})$ but not over $\mathbb{Q}$.

We recall that $E$ is modular and call the attached eigenform $f_{E}$ with $%
q $-expansion
\begin{equation*}
f_{E}=q+\sum\limits_{n=2}^{\infty }a_{n}q^{n}.
\end{equation*}
This is a newform in $S_{2}(N_{E},\chi _{1})$ where
$S_{2}(N_{E},\chi _{1})$ is the space of cusp forms of weight $2$,
level $N_{E}$ and $\chi _{1}$ is the trivial character.

The attached eigenform of $E_{D}$ is the twist of $f_{E}$ by the
quadratic character $\chi _{D}$ $:$ $f_{E_{D}}:=f_{E}\otimes \chi
_{D}=\sum\limits_{n=1}^{\infty }\chi _{D}(n)a_{n}q^{n}\in
S_{2}(N_{E_{D}})$ (and $N_{E_{D}}$ divides $N_{E}.D^{2}$). So the
Hasse-Weil $L-$function of $E_{D}$ is
\begin{equation*}
L_{E_{D}}(s)=\sum\nolimits_{n=1}^{\infty }\chi _{D}(n)a_{n}n^{-s}.
\end{equation*}

In this paper, we shall give numerical data related to the order
of the Selmer groups of twist families $\{E_{D}\}.$ In particular
we are interested in the number of twists for which there are
infinitely many points in $E_{D}(\mathbb{Q})$. Recall the theorem
of Mordell which states that
\begin{equation*}
E(\mathbb{Q})\cong E(\mathbb{Q})_{tor}\times\mathbb{Z}^{r},
\end{equation*}
where the torsion subgroup $E(\mathbb{Q})_{tor}$ is finite and the
rank $r$ of $E(\mathbb{Q})$ is a non-negative integer.

For any given elliptic curve, it is possible to describe quite
precisely the torsion subgroup \cite{Silverman}. The rank is much
more difficult to compute, and in general there is no known
procedure which is guaranteed to yield an answer. But if the rank
of $E$ is positive then a celebrated theorem of Kolyvagin states
that $L_{E}(1)=0$. So we are sure that if $L_{E}(1)\neq 0$ then
$E(\mathbb{Q})$ is finite. The converse result is not known today
but it should be true. One part of the celebrated \emph{Birch and
Swinnerton-Dyer Conjecture (BSD)} is that the order of vanishing
of $L_{E}(s)$ at $s=1$ ("the analytic rank") is equal to the rank
of $E(\mathbb{Q}).$ BSD states much more. It interprets the value
of the first non-vanishing
derivative of $L_{E}$ at $s=1$ in terms of arithmetical objects attached to $%
E.$ We shall be interested in this prediction only in the case
that the analytic rank of $E$ is $0$.

We begin defining the Selmer and the Tate-Shafarevich group of
elliptic curves by using the Kummer sequence of elliptic curves:
Let $E$ be an elliptic curve over $\mathbb{Q}$. Let
$\overline{\mathbb{Q}}$ be an algebraic closure of $\mathbb{Q}$
and $G_{\mathbb{Q}}:=Aut_{\mathbb{Q}}(\overline{\mathbb{Q}})$ the
absolute Galois group of $\mathbb{Q}$. Consider the abelian group
$E(\overline{\mathbb{Q}})$ of all points on $E$ defined over
$\overline{\mathbb{Q}}$. One can consider the Galois cohomology
groups $H^{m}(G_{\mathbb{Q}},E(\overline{\mathbb{Q}}))$ for $m$
$\in \mathbb{N}$.

For all $n\in\mathbb{N}$, we have the exact sequence of
$G_{\mathbb{Q}}-$modules
\begin{equation*}
0\longrightarrow E(\overline{\mathbb{Q}})[n]\longrightarrow
E(\overline{\mathbb{Q}})\overset{n}{\longrightarrow
}E(\overline{\mathbb{Q}})\longrightarrow 0.
\end{equation*}

As it is well known \cite{Serre}, there is an associated long
exact sequence of Galois cohomology groups. We need a consequence
of the beginning of this sequence \cite{Silverman}
\begin{equation*}
0\rightarrow E(\mathbb{Q})/nE(\mathbb{Q})\rightarrow
H^{1}(G_{\mathbb{Q}},E(\overline{\mathbb{Q}})[n])\overset{\alpha
}{\rightarrow
}H^{1}(G_{\mathbb{Q}},E(\overline{\mathbb{Q}}))[n]\rightarrow 0.
\end{equation*}
This sequence is called the \emph{Kummer Sequence} associated to
$E$. For each prime $p$ we choose an extension of the
corresponding $p-$adic valuation. Let $G_{p}$ be the corresponding
decomposition group in $G_{\mathbb{Q}}$ which is in a canonical
way isomorphic to $G_{\mathbb{Q}_{p}}$. Let $\gamma _{p,n}$ be the
restriction map from
$H^{1}(G_{\mathbb{Q}},E(\overline{\mathbb{Q}})[n])$ to
$H^{1}(G_{p},E(\overline{\mathbb{Q}}_{p}))[n]$ and $P$ the set of
primes. The \emph{Tate-Shafarevich group} of $E$ is denoted by
$Sha_{\mathbb{Q}}(E)$ and defined by
\begin{equation*}
Sha_{\mathbb{Q}}(E):=\underset{n\in \mathbb{N}}{\bigcup
}Sha_{\mathbb{Q}}(E)[n],
\end{equation*}
where
\begin{equation*}
Sha(E)[n]:=\bigcap\limits_{p\in P}\ker (\gamma _{p,n}).
\end{equation*}
The \emph{Selmer group} of $E$ is denoted by $S_{\mathbb{Q}}(E)$
and defined by
\begin{equation*}
S_{\mathbb{Q}}(E):=\underset{n\in\mathbb{N}}{\bigcup
}S_{\mathbb{Q}}(E)[n],
\end{equation*}
where
\begin{equation*}
S_{\mathbb{Q}}(E)[n]:=\alpha ^{-1}(Sha_{\mathbb{Q},S}(E)[n]).
\end{equation*}
So we have the exact sequence
\begin{equation*}
0\longrightarrow E(\mathbb{Q})/nE(\mathbb{Q})\longrightarrow
S_{\mathbb{Q}}(E)[n]\longrightarrow
Sha_{\mathbb{Q}}(E)[n]\longrightarrow 0.
\end{equation*}

We are now ready to state the part of BSD which is of importance
for us.
\begin{conjecture}
\cite{BSD} $L_{E}(1)\neq 0$ iff $E(\mathbb{Q})$ is finite group,
and then the Selmer group of $E$ is finite and the following
equality holds:
\begin{equation*}
L_{E}(1)=\left( \underset{E^{0}(\mathbb{R})}{\int }|\omega
_{E}|\right) \left( \underset{p|N.\infty }{\prod }c_{p}\right)
\frac{\#S_{\mathbb{Q}}(E)}{\#(E(\mathbb{Q}))^{3}},
\end{equation*}
where $E^{0}(\mathbb{R})$ is the connected component of
$E(\mathbb{R})$, $\omega _{E}$ is the N\'{e}ron differential of
$E$, $c_{\infty }=[E(\mathbb{R}):E^{0}(\mathbb{R})]$, and for
primes $p$, $c_{p}=[E(\mathbb{Q}_{p}):E^{0}(\mathbb{Q}_{p})]$. The
numbers $c_{p}$ are called \emph{local} \emph{Tamagawa numbers}.
\end{conjecture}

We remark that all terms different from the order of the Selmer
group are computable more or less easily. But in some special
cases it is possible to compute the order of the Selmer groups
(sometimes one has to assume its finiteness), and then one can
verify BSD. So there is numerical evidence for its truth.

\textbf{Convention: }Without further notice we always shall
\textbf{assume } in this paper that BSD holds and use the analytic
theory of modular forms to compute both the order of
$S_{\mathbb{Q}}(E)$ and $L_{E}(1)$ conditionally.

A good test for the exactness of algorithms is a result of Cassels
for the order of $S_{\mathbb{Q}}(E)$:

\begin{theorem}
\cite{Cassels} Let $E/\mathbb{Q}$ be an elliptic curve. There
exists an alternating, bilinear pairing
\begin{equation*}
\Gamma :Sha_{\mathbb{Q}}(E)\times
Sha_{\mathbb{Q}}(E)\longrightarrow \mathbb{Q}/\mathbb{Z}\end{equation*}%
whose kernel is precisely the group of divisible elements of
$Sha$.

In particular if $S_{\mathbb{Q}}(E)$ is finite, then
$k=\#S_{\mathbb{Q}}(E)/\#(E(\mathbb{Q}))$ is a perfect square.
\end{theorem}

\section{Waldspurger Theorem and Its Consequences}

\subsection{Statement of Waldspurger's Theorem}

Assume that the rank of $E$ is equal to zero. As said above one
can compute the order of the Selmer group and hence of the
Tate-Shafarevich group of $E$ by using BSD. Note that the local
Tamagawa numbers $c_{p}$ as well as $\omega _{E}$ can be computed
easily (the latter value is transcendental and hence has to be
computed up to a desired precision). The most time consuming item
is the computation of $L_{E}(1)$. For this, there is a routine in
the computer algebra system MAGMA \cite{Magma}.

It turns out that computing $L_{E}(1)$ with the necessary
precision (again this is a transcendental) for an elliptic curve
$E$ with large conductor takes a long time.

For instance, computing $L_{E}(1)$ for the elliptic curve%
\begin{equation*}
E:y^{2}=x^{3}-87662765543106x+572205501116432432042932656
\end{equation*}%
which has conductor $11520793560025904$, one needs at least $1000$
hours in
a laptop computer$\footnote{%
With the properties: Intel Core 2 Duo Mobile, 2GB DDR2, 2.00GHz}$
with MAGMA which doesn't guarantee to answer. Another hard
numerical problem is to decide by computation whether
$L_{E}(1)=0$.

The situation is much better in families of twists of a given
elliptic curve. The elliptic curve $E$ from above is a member of
such family, and we shall see in Section $3.3$ how this can be
used to accelerate the computation dramatically. The reason is
\emph{Waldspurger's Theorem }which is crucial for our work:

\begin{theorem}
\cite{Waldspurger} Let $E$ be an elliptic curve over $\mathbb{Q}$
with attached new form $f_{E}$. Assume that $F_{E}\in
S_{3/2}(N^{\prime },\chi _{1})$ is an eigenform and
\textbf{S}$(F_{E})=f_{E}$ where \textbf{S} is the Shimura-Shintani
lifting.

Let $a_{n}$ be the $n-$th Fourier coefficient of $F_{E}$. Then for
square-free natural numbers $n$ and $n_{0}$ with $n\equiv n_{0}$ mod $%
\prod\limits_{p|N^{\prime }}\mathbb{Q}_{p}^{\ast ^{2}}$ and
$n.n_{0}$ prime to $N^{\prime }$ we have
\begin{equation*}
a_{n_{0}}^{2}\sqrt{n}L_{E_{-n}}(1)=a_{n}^{2}\sqrt{n_{0}}L_{E_{-n_{0}}}(1),
\end{equation*}
Hence we get: If $a_{n_{0}}\neq 0$ then $L_{E_{-n}}(1)$ is
determined by $a_{n}$, $a_{n_{0}}$ and $L_{E-n_{0}}(1).$

In particular, $L_{E_{-n}}(1)=0$ for all $n\equiv n_{0}$ mod $%
\prod\limits_{p|N^{\prime }}\mathbb{Q}_{p}^{\ast ^{2}}$ iff
$L_{E-n_{0}}(1)=0$, and else $L_{E_{-n}}(1)\neq 0$ iff $a_{n}\neq
0$.
\end{theorem}

\begin{corollary}
\cite{ABF} Assume that the Birch and Swinnerton-Dyer conjecture holds for $%
E_{-n}$ and $E_{-n_{0}}$, $n$, $n_{0}$ as in the theorem and that
$a_{n_{0}}\cdot L_{E-n_{0}}(1)\neq 0$. Then $E_{-n}(\mathbb{Q})$
is finite iff $a_{n}\neq 0$ and
\begin{equation*}
\#S_{\mathbb{Q}}(E_{-n})=d(n,n_{0})\cdot
\#S_{\mathbb{Q}}(E_{-n_{0}})\frac{a_{n}^{2}}{a_{n_{0}}^{2}},
\end{equation*}
where $d(n,n_{0})$ is easily computed as explained in Subsection
$3$. $2$
and essentially a power of $2$ depending on the divisor structure of $n$, $%
n_{0}$.
\end{corollary}

\subsection{Computing $d(n,n_{0})$}

We continue to assume that $E_{-n}$ and $E_{-n_{0}}$ are twists of
$E$ satisfying the conditions of Corollary $4$.

We want to compute the numbers $d(n,n_{0})$. By definition
$d(n,n_{0})$ depends on the Tamagawa numbers and the torsion
subgroups of the two elliptic curves $E_{-n}$ and $E_{-n_{0}}$.

To be explicit one has to use some easy facts about twists of
elliptic curves. $d(n,n_{0})$ does not depend on the real period
since the twisting factors $n_{0}$ and $n$ are odd and congruent
modulo $4$ and $\omega _{E-n}/\omega
_{E-n_{0}}=\sqrt{n_{0}}/\sqrt{n}$ and hence cancel in the formula
in Corollary $4$.

\textbf{Independence of torsion elements:}

It is well known and obvious that for all pairs of elliptic curve
$E$ and twists $E_{-n}$ we have
$E(\mathbb{Q})[2]=E_{-n}(\mathbb{Q})[2]$.

Moreover for given $E$ there are only finitely many (in fact only
very few) twists of $E$ which have torsion points of order $>2$
over $\mathbb{Q}$. Avoiding these twists is easy and so one can
assume without falsifying the statistic, that all members of the
twist families have only $\mathbb{Q}-$rational torsion points of
order dividing $2.$ In fact, in the chosen examples below this
holds for all non-trivial twists of the treated curves $E $. So we
can assume that the order of $E_{-n_{0}}(\mathbb{Q})$ is equal to
the order of $E_{-n}(\mathbb{Q})$ and hence $d(n,n_{0})$ is
independent of torsion elements.

The next observation is that the groups of connected components of
twists of an elliptic curve $E$ over the reals are equal, and so
$d(n,n_{0})$ is computed by looking at the non-Archimedean
Tamagawa numbers.

Let us denote the Tamagawa numbers for $E_{-n}$ at a prime $p$ by
$c_{n,p}$ and the Tamagawa numbers for $E_{-n_{0}}$ at a prime $p$
by $c_{n_{0},p}.$

First observe that for $p$ prime to $n_{0}\cdot n\cdot N^{\prime
}$ both twists have good reduction modulo $p$ and so the Tamagawa
numbers are equal to $1.$

By assumption $-n_{0}$ and $-n$ lie in the same class of squares
in all
completions with respect to divisors of $N_{E}^{\prime }$ and so the N\'{e}%
ron models are equal at all primes dividing $N_{E}^{\prime }$.

Now let $p$ be a divisor of, say, $n$ prime to $N_{E}^{\prime }$.
Since $E$ has good reduction modulo $p$ and $p$ is odd we can use
the table $15.1$ in
\cite{Silverman}, p. $359$ to see that $E$ has Kodaira symbol $\mathcal{I}%
_{0}$ and so $c_{n,p}=4$. The same result holds of course for
prime divisors of $n_{0}.$

Hence we get

\begin{lemma}
Let $E$ be an elliptic curve and $E_{-n}$ and $E_{-n_{0}}$ be
twists of $E$ with
$\#E_{-n}(\mathbb{Q})=\#E_{-n_{0}}(\mathbb{Q})<\infty $ and
$n.n_{0}$ prime to $N_{E}^{\prime }$ and a square in all
completions with respect to divisors of $N_{E}^{\prime }.$ Then
\begin{equation*}
d(n,n_{0})=\frac{\prod c_{n_{0},p}}{\prod
c_{n,p}}=\frac{4^{\#div(n_{0})}}{4^{\#div(n)}},
\end{equation*}
where $\#div(-)$ denotes the number of prime divisors of $-$.
\end{lemma}

For using Waldspurger Theorem for members of the twist family
$\{E_{-n}\}$ one has to find an eigenform $F_{E}$ as above. Then
one has to implement a fast algorithm for computing the Fourier
coefficients of $F_{E}$ in a large range.

\subsection{Computing Fourier Coefficients and the Selmer
Group}

Recall the situation. We have an elliptic curve $E$ with eigenform
$f_{E}$ and the Shimura-Shintani lift $F_{E}$ given in a concrete
way. In particular we shall consider the following examples from
\cite{Frey}:
\begin{equation*}
\begin{tabular}{|c|c|}
\hline $F_{E}$ & $E$ \\ \hline \multicolumn{1}{|l|}{$(\Theta
(X^{2}+11Y^{2})-\Theta
(3X^{2}+2XY+4Y^{2})).\Theta _{id,11}$} & \multicolumn{1}{|l|}{$11a1$} \\
\hline \multicolumn{1}{|l|}{$(\Theta (X^{2}+14Y^{2})-\Theta
(2X^{2}+7Y^{2})).\Theta _{id,14}$} & \multicolumn{1}{|l|}{$14a1$}
\\ \hline \multicolumn{1}{|l|}{$(\Theta
(3X^{2}-2XY+23Y^{2})-\Theta
(7X^{2}+6XY+11Y^{2})).\Theta _{id,17}$} & \multicolumn{1}{|l|}{$17a1$} \\
\hline \multicolumn{1}{|l|}{$(\Theta (X^{2}+20Y^{2})-\Theta
(4X^{2}+5Y^{2})).\Theta _{id,20}$} & \multicolumn{1}{|l|}{$20a1$}
\\ \hline \multicolumn{1}{|l|}{$(\Theta (X^{2}+17Y^{2})-\Theta
(2X^{2}+2XY+9Y^{2})).\Theta _{id,17}$} & \multicolumn{1}{|l|}{$34a1$} \\
\hline
\end{tabular}%
\end{equation*}%
where $\Theta (.)$ is the theta series of a binary quadratic form and $%
\Theta _{\psi ,t}:=\underset{n=-\infty }{\overset{\infty
}{\sum}}\psi (n)q^{tn^{2}}$ is a Fourier series for the Dirichlet
character $\psi .$

The elliptic curve $E$ is given as in Cremona's Table
\cite{Cremona}.

Let $F_{E}\in S_{3/2}(N^{\prime },\chi _{1})$ as above with
Fourier expansion $\underset{n=1}{\overset{\infty }{\sum
}}a_{n}q^{n}$.

\textbf{Strategy}

1) Calculate the $q-$expansion of $F_{E}$ up to an upper bound
$M$, construct the list $L:=\{(n,a_{n})|n\in \{1,\cdots,M\}$
squarefree$\}$.

2) Choice of Congruence Classes: To apply Waldspurger's theorem we
compare twists with twist factors $-n$, $-n_{0}$ with $n$ and
$n_{0}$ odd and prime to $N^{\prime }$ which are congruent modulo
$\prod\limits_{p|N^{\prime }}\mathbb{Q}_{p}^{\ast ^{2}}$. This is
satisfied if $n\equiv n_{0}$ $mod$ $8\cdot \underset{2\neq
p|N_{E}}{\prod }p$ and hence we shall investigate twist families
with twist factors in such congruence classes. First we determine
the twist families (with respect to the above congruences) which
consist of \emph{odd} elliptic curves and so have positive
analytic rank by looking at the parity of the twist characters. We
delete these congruence classes.

We simplify the situation in the cases $N_{E}=11$ and $N_{E}=17$.
To apply Waldspurger's theorem we have to look at congruence
classes modulo $88$ and respectively $136$. We check that for all
pairs of these congruence classes which become equal modulo $44$
respectively $68$ there are $n_{0}n_{0}^{\prime }$ with the same
number of prime divisors, the Fourier coefficients
$a_{n_{0}}=a_{n_{0}^{\prime }}$ and the same order of the Selmer
groups and hence we can investigate in these cases twist families
with families with twist factors running over congruences modulo
$44$ respectively $68$.

We list the resulting congruence classes in Table A.

\begin{center}
\begin{tabular}{|c|c|c|}
\hline Elliptic Curve & Modulo & $n_{0}$ \\ \hline $11a1$ & $44$ &
$1,3,5,15,23,31,37$ \\ \hline $14a1$ & $56$ &
$1,15,23,29,37,39,53$ \\ \hline $17a1$ & $68$ & $3,7,11,23,31,39$
\\ \hline $20a1$ & $40$ & $1,21,29$ \\ \hline
$34a1$ & $136$ & $1,13,19,21,33,35,43,53,59,67,69,77,$ \\
&  & $83,89,93,101,115,117,123$ \\ \hline
\end{tabular}

Table A.
\end{center}

3) For the integer $M$ and fixed $n_{0}$ calculate
\begin{equation*}
x_{n_{0}}(M):=\#\{n:n\leq M,\ n\text{ is square-free,}\ n\equiv
n_{0} \ (mod N^{\prime })\},
\end{equation*}%
\begin{equation*}
s_{n_{0},0,E}(M):=\#\{n:n\leq M\text{, }n\text{ is square-free,
}n\equiv n_{0}\text{ }(mod\text{ }N^{\prime })\text{,
}a_{n}=0\}\text{,}
\end{equation*}%
and plot the function $s_{n_{0},0,E}(M)/x_{n_{0}}(M)$.

4) For $n_{0}$, find $\alpha \in\mathbb{R}$ and $\epsilon \in
\lbrack -0.02,0.02]$ such that

\begin{equation*}
\sigma (x_{n_{0}}(M)):=\alpha \frac{(\log \log (x_{n_{0}}(M)))^{1+\epsilon }%
}{\log (x_{n_{0}}(M))}.
\end{equation*}%
approximates $s_{n_{0},0,E}(M)/x_{n_{0}}(M)$ "well".

5) If $a_{n_{0}}=0$ then replace $n_{0}$ by the minimal $n$ in the
congruence class such that $a_{n_{0}}\neq 0$. Calculate
$L_{E_{-n_{0}}}(1)$, $\#E_{-n_{0}}(\mathbb{Q})_{tors}$ and
$\#S_{\mathbb{Q}}(E_{-n_{0}})$ by using the BSD-conjecture and
$f_{E}.$

6) For $n_{0}$ and $n\leq M$, compute $d(n,n_{0})$ as described in
Subsection $3.2$.

7) For $n_{0}$ and $n\leq M$, compute
\begin{equation*}
s_{E_{-n}}:=\frac{\#S_{\mathbb{Q}}(E_{-n_{0}})\cdot a_{n}^{2}\cdot
d(n,n_{0})}{a_{n_{0}}^{2}}.
\end{equation*}%
which is, according to the BSD-conjecture, the order of
$S_{\mathbb{Q}}(E_{-n})$.

8) Compute $t:=\#E(\mathbb{Q})$. Recall that twisting $E$ doesn't
change the order of the torsion subgroup of $E$.

9) For $M,k,t$ and $n_{0}$ compute

\begin{equation*}
s_{n_{0},k,E}(M):=\#\{n:n\leq M\text{, }n\text{ is square-free,
}n\equiv
n_{0}\text{ }(mod\text{ }N^{\prime })\text{, }\frac{s_{E_{-n}}}{t}=k\}%
\text{.}
\end{equation*}

10) For $n_{0}$, plot the function
$s_{n_{0},k,E}(M)/x_{n_{0}}(M)$.

11) For $n_{0}$, find $\alpha \in\mathbb{R}$ and $\epsilon \in
\lbrack -0.02,0.02]$ such that
\begin{equation*}
\sigma (x_{n_{0}}(M))=\alpha \frac{(\log \log (x_{n_{0}}(M)))^{1+\epsilon }}{%
\log (x_{n_{0}}(M))}
\end{equation*}%
approximates $s_{n_{0},k,E}(M)/x_{n_{0}}(M)$ "well".

\begin{remark}
\emph{1)} All data can be found in \emph{http://homepage.uludag.edu.tr/\symbol{126}%
inam/}

\emph{2)} Having computed $\#S_{\mathbb{Q}}(E_{-n})$, $d(n,n_{0})$
and $L_{E_{-n_{0}}}(1)$, one can use the BSD-conjecture again to
compute $L_{E_{-n}}(1)$ as
\begin{equation*}
L_{E_{-n}}(1)=\frac{L_{E_{-n_{0}}}(1)\cdot
\#S_{\mathbb{Q}}(E_{-n})}{\#S_{\mathbb{Q}}(E_{-n_{0}})\cdot
d(n,n_{0})}.
\end{equation*}%
This is much faster than to compute $L_{E_{-n}}(1)$ directly. We
have included these values in our lists.
\end{remark}

We come back to the example in Section $3.1$. Recall that we
wanted to compute the value of $L_{E}(1)$ of the curve
\begin{equation*}
E:y^{2}=x^{3}-87662765543106x+572205501116432432042932656.
\end{equation*}
It is the twist of the elliptic curve $11a1$ with the twist factor $%
n=8090677 $, and by the method described above, we get very fast%
\begin{equation*}
L_{E}(1)=2.10072023061090418110927626775
\end{equation*}
approximately in $360$ seconds.

We now fix an elliptic curve $E$ as well as $n_{0}$ and $k$.

We sketch how to determine an approximation function for
$q_{n_{0},k,E}$. We choose $\alpha $ and $\varepsilon $ in the
following way: In this work,
using the data obtained up to the bound $M=10^{7}$, we construct a family $%
\{I_{i}\}$ of subintervals of $I:=[0,M]$ defined by
$I_{i}=[0,50000i]$ for $i=1,2,\cdots,200$ such that
\begin{equation*}
I_{1}\subseteq I_{2}\subseteq \cdots \subseteq I_{200}=I.
\end{equation*}

First we calculate the value
$\frac{s_{n_{0},k,E}(M_{i})}{x_{n_{0}}(M_{i})}$ and afterwards
$\frac{\log (\log (x_{n_{0}}(M_{i}))}{x_{n_{0}}(M_{i})}$.
Comparing these two values, the constant $\alpha _{i}$ can be
obtained for each $I_{i}$. Using the weighted average for the
constants $\alpha _{i}$ we
determine $\alpha $ (depending on $M$). By means of these constants $\alpha $%
, we compare the values
\begin{equation*}
\frac{s_{n_{0},k,E}(M_{i})}{x_{n_{0}}(M_{i})}\text{ and }\alpha
\frac{\log (\log (x_{n_{0}}(M_{i}))}{x_{n_{0}}(M_{i})}.
\end{equation*}%

After this step we choose for fine-tuning $\varepsilon \in $ $%
[-0.02,0.02]$ such that the approximation is getting better.

\section{Numerical Results}

We use the considerations of Subsection $3.3$ for extensive
computations and observe that in all our examples the functions
\begin{equation*}
q_{n_{0},k,E}:=\frac{s_{n_{0},k,E}}{x_{n_{0}}}(M)
\end{equation*}%
are fairly well approximated by
\begin{equation*}
\alpha \frac{\log (\log (x_{n_{0}}(M))^{1+\varepsilon }}{\log
(x_{n_{0}}(M))}
\end{equation*}%
where $\alpha >0$ and $\varepsilon \in \lbrack -0.02,0.02]$.

This observation confirms predictions stated by Birch and lead to
\begin{conjecture}
For all elliptic curves $E$, $E^{\prime }$ over $\mathbb{Q}$, all
$n_{0}$, $n_{0}^{\prime }$ satisfying the conditions of Theorem
$3$ and all $k$, $k^{\prime }$ the asymptotic behavior of
$\frac{q_{n_{0},k,E}}{q_{n_{0}^{\prime },k^{\prime },E^{\prime
}}}$ is well approximated by a constant times a factor $\log (\log
(x(M))^{\delta }$ where $x(M)$ is the number of square-free
numbers $\leq M$ and $\delta $ is a real number with small
absolute value.
\end{conjecture}

Of course one should be much more precise and predict how the
factor depends on the parameters. In our context we shall restrict
ourselves to a discussion of the reals $\alpha $ we get out of our
data by the approximation process in the algorithm described
above.

\subsection{Observations}

Conjecture $7$ predicts that the type of the approximation
function is independent of $k$. But the constants $\alpha $ vary
and so a finer analysis seems necessary in order to find reasons
or patterns for the size of $\alpha .$

But let us begin with a word of caution. {\normalsize In our
examples we computed }$\alpha ${\normalsize \ for }$k\leq 961$.
For large $k$ we do not have enough material for any statistical
statement. (The record, $k=68121$ occurs just one time).

We now discuss examples of the weighted average values which are
given in Section $4$.$3$.

\subsubsection{Some Examples}

Considering the values of $\alpha $ in our examples we see that it
depends significantly on the congruence classes modulo $4$
respectively $8$.

\emph{Case }$1.$ For the elliptic curve $11a1$, we have the class $K$ :$%
=\{1,5,37\}$ which have members congruent to $1$ modulo $4$ and the class $%
L:=\{3,15,23,31\}$ with members congruent to $3$ modulo $4$. As
example, we take $k=1$ and see that the values of $\alpha $ for
$n_{0}$ in class $K$ are respectively $0.296$, $0.299$, $0.300$,
where as for $n_{0}$ in class $L$
they are respectively $0.458$, $0.459$, $0.469$, and $0.464$ (See Section $4$%
.$3$).

\emph{Case }$2.$ For the elliptic curve $14a1$, we have one class
modulo $8$ in which all twists are odd curves, namely $n$
congruent to $3$ modulo $8$. The other classes modulo $8$ separate
our congruence classes modulo $56$
into $K:=\{1\}$, $L:=\{29,37,53\}$, $M:=\{15,23,39\}$. As example we take $%
k=9$ and get the values of $\alpha $ for $n_{0}$ in class $K$ are
$0.485$, for $n_{0}$ in class $L$ are respectively $0.195$,
$0.185$ and $0.192$
whereas for $n_{0}$ in class $M$ they are respectively $0.559$, $0.568$ and $%
0.536.$

\emph{Case }$3.$ For the elliptic curve $17a1$. In this case all
congruence classes modulo $68$ which are congruent to $1$ modulo
$4$ are odd, and for all classes congruent $3$ modulo $4$ the
values $\alpha $ are around $0.33$ and hence of the same size for
$k=0$ as an example.

\emph{Case }$4.$ For the elliptic curve $20a1$, all congruence
classes modulo $40$ which are congruent to $3$ modulo $4$ are odd,
and for all classes congruent $1$ modulo $4$ the values $\alpha $
are around $0.1$ and hence of the same size for $k=225$ as an
example.

\emph{Case }$5$. For the elliptic curve $34a1$, we looked at
congruence classes modulo $8$. If $n$ is congruent $7$ modulo $8$
we get odd curves.
The other classes modulo $8$ consist of $K:=\{1,33,89\}$, $%
L:=\{19,35,43,59,67,83,115,123\}$,
$M:=\{21,53,69,77,93,101,117\}$. Again
take $k=1$, then the values of $\alpha $ are around $0.38$ and for every $%
n_{0}\in L$, the values of $\alpha $ are around $0.47$ and for every $%
n_{0}\in M$, the values of $\alpha $ are around $0.41.$

\subsection{Some Actual Values}

We give some examples of actual values. Here, all the notation
given before is valid and $\sigma (x_{n_{0}}(M))$ is defined by
$\sigma (x_{n_{0}}(M)):=\alpha \frac{(\log \log
(x_{n_{0}}(M)))^{1+\epsilon }}{\log (x_{n_{0}}(M))}$. The values
of $\alpha $ are given in Subsection $4.3$.

For the elliptic curve $11a1$, $n_{0}=3,k=4$ and $\varepsilon
=0.005$ we have

\begin{center}
\begin{tabular}{|c|l|l|}
\hline
$M$ & $s_{n_{0},k}(x_{n_{0}}(M))/x_{n_{0}}(M)$ & $\sigma (x_{n_{0}}(M))$ \\
\hline \multicolumn{1}{|l|}{$50000$} & $0.106452$ & $0.099558$ \\
\hline \multicolumn{1}{|l|}{$1500000$} & $0.074267$ & $0.066593$
\\ \hline \multicolumn{1}{|l|}{$3000000$} & $0.066195$ &
$0.062381$ \\ \hline \multicolumn{1}{|l|}{$4000000$} & $0.062997$
& $0.060786$ \\ \hline \multicolumn{1}{|l|}{$5000000$} &
$0.060743$ & $0.059604$ \\ \hline \multicolumn{1}{|l|}{$10000000$}
& $0.053981$ & $0.056209$ \\ \hline
\end{tabular}
\end{center}

For the elliptic curve $14a1$, $n_{0}=1,k=16$ and $\varepsilon
=0.005$ we have
\begin{center}
\begin{tabular}{|c|l|l|}
\hline
$M$ & $s_{n_{0},k}(x_{n_{0}}(M))/x_{n_{0}}(M)$ & $\sigma (x_{n_{0}}(M))$ \\
\hline \multicolumn{1}{|l|}{$100000$} & $0.082313$ & $0.09115$ \\
\hline \multicolumn{1}{|l|}{$1400000$} & $0.066638$ & $0.066978$
\\ \hline \multicolumn{1}{|l|}{$2000000$} & $0.06462$ & $0.064665$
\\ \hline \multicolumn{1}{|l|}{$5000000$} & $0.060241$ &
$0.059394$ \\ \hline \multicolumn{1}{|l|}{$8000000$} & $0.056955$
& $0.057011$ \\ \hline \multicolumn{1}{|l|}{$10000000$} &
$0.055412$ & $0.055946$ \\ \hline
\end{tabular}
\end{center}

For the elliptic curve $17a1$, $n_{0}=7,k=324$ and $\varepsilon
=0.005$ we have
\begin{center}
\begin{tabular}{|c|l|l|}
\hline
$M$ & $s_{n_{0},k}(x_{n_{0}}(M))/x_{n_{0}}(M)$ & $\sigma (x_{n_{0}}(M))$ \\
\hline \multicolumn{1}{|l|}{$100000$} & $0$ & $0.016898$ \\ \hline
\multicolumn{1}{|l|}{$5000000$} & $0.009965$ & $0.010903$ \\
\hline \multicolumn{1}{|l|}{$6000000$} & $0.010771$ & $0.010726$
\\ \hline \multicolumn{1}{|l|}{$7000000$} & $0.011213$ & $0.01058$
\\ \hline \multicolumn{1}{|l|}{$8000000$} & $0.011651$ &
$0.010458$ \\ \hline \multicolumn{1}{|l|}{$10000000$} & $0.012272$
& $0.010259$ \\ \hline
\end{tabular}
\end{center}
For the elliptic curve $20a1$, $n_{0}=1,k=100$ and $\varepsilon
=0.005$ we have

\begin{center}
\begin{tabular}{|c|l|l|}
\hline
$M$ & $s_{n_{0},k}(x_{n_{0}}(M))/x_{n_{0}}(M)$ & $\sigma (x_{n_{0}}(M))$ \\
\hline \multicolumn{1}{|l|}{$500000$} & $0.026748$ & $0.038045$ \\
\hline \multicolumn{1}{|l|}{$3000000$} & $0.029427$ & $0.031896$
\\ \hline \multicolumn{1}{|l|}{$5000000$} & $0.029764$ &
$0.030491$ \\ \hline \multicolumn{1}{|l|}{$6000000$} & $0.029958$
& $0.030019$ \\ \hline \multicolumn{1}{|l|}{$7000000$} &
$0.030039$ & $0.029632$ \\ \hline \multicolumn{1}{|l|}{$10000000$}
& $0.030132$ & $0.028772$ \\ \hline
\end{tabular}
\end{center}

For the elliptic curve $34a1$, $n_{0}=1,k=36$ and $\varepsilon
=0.005$ we have
\begin{center}
\begin{tabular}{|c|l|l|}
\hline
$M$ & $s_{n_{0},k}(x_{n_{0}}(M))/x_{n_{0}}(M)$ & $\sigma (x_{n_{0}}(M))$ \\
\hline \multicolumn{1}{|l|}{$3000000$} & $0.066667$ & $0.071691$
\\ \hline \multicolumn{1}{|l|}{$5000000$} & $0.069827$ &
$0.069099$ \\ \hline \multicolumn{1}{|l|}{$6000000$} & $0.068564$
& $0.067903$ \\ \hline \multicolumn{1}{|l|}{$7000000$} & $0.0682$
& $0.066922$ \\ \hline \multicolumn{1}{|l|}{$8000000$} &
$0.067812$ & $0.066096$ \\ \hline \multicolumn{1}{|l|}{$10000000$}
& $0.067339$ & $0.064759$ \\ \hline
\end{tabular}%
\end{center}

\subsection{ Weighted average values of $\protect\alpha $'s}

\begin{center}
\begin{tabular}{|c|c|c|c|c|c|c|}
\hline $E$ & $n_{0}$ & $0$ & $1$ & $4$ & $9$ & $16$ \\ \hline
\multicolumn{1}{|l|}{$11a1$} & \multicolumn{1}{|l|}{$1$} &
\multicolumn{1}{|l|}{$0.140221$} &
\multicolumn{1}{|l|}{$0.295669$} &
\multicolumn{1}{|l|}{$0.204751$} &
\multicolumn{1}{|l|}{$0.309679$} &
\multicolumn{1}{|l|}{$0.184184$} \\ \hline
\multicolumn{1}{|l|}{$11a1$} & \multicolumn{1}{|l|}{$3$} &
\multicolumn{1}{|l|}{$0.214424$} &
\multicolumn{1}{|l|}{$0.458141$} &
\multicolumn{1}{|l|}{$0.296646$} &
\multicolumn{1}{|l|}{$0.445157$} &
\multicolumn{1}{|l|}{$0.244439$} \\ \hline
\multicolumn{1}{|l|}{$11a1$} & \multicolumn{1}{|l|}{$5$} &
\multicolumn{1}{|l|}{$0.139959$} &
\multicolumn{1}{|l|}{$0.299438$} &
\multicolumn{1}{|l|}{$0.199569$} &
\multicolumn{1}{|l|}{$0.308824$} &
\multicolumn{1}{|l|}{$0.186938$} \\ \hline
\multicolumn{1}{|l|}{$11a1$} & \multicolumn{1}{|l|}{$15$} &
\multicolumn{1}{|l|}{$0.211029$} &
\multicolumn{1}{|l|}{$0.458734$} &
\multicolumn{1}{|l|}{$0.299005$} &
\multicolumn{1}{|l|}{$0.442075$} &
\multicolumn{1}{|l|}{$0.244968$} \\ \hline
\multicolumn{1}{|l|}{$11a1$} & \multicolumn{1}{|l|}{$23$} &
\multicolumn{1}{|l|}{$0.208441$} &
\multicolumn{1}{|l|}{$0.468673$} &
\multicolumn{1}{|l|}{$0.304083$} &
\multicolumn{1}{|l|}{$0.440648$} & \multicolumn{1}{|l|}{$0.23676$}
\\ \hline \multicolumn{1}{|l|}{$11a1$} &
\multicolumn{1}{|l|}{$31$} & \multicolumn{1}{|l|}{$0.208064$} &
\multicolumn{1}{|l|}{$0.46357$} & \multicolumn{1}{|l|}{$0.205327$}
& \multicolumn{1}{|l|}{$0.441027$} &
\multicolumn{1}{|l|}{$0.23975$} \\ \hline
\multicolumn{1}{|l|}{$11a1$} & \multicolumn{1}{|l|}{$37$} &
\multicolumn{1}{|l|}{$0.14234$} & \multicolumn{1}{|l|}{$0.300449$}
& \multicolumn{1}{|l|}{$0.205431$} &
\multicolumn{1}{|l|}{$0.314227$} & \multicolumn{1}{|l|}{$0.18237$}
\\ \hline
\end{tabular}

\begin{tabular}{|c|c|c|c|c|c|c|}
\hline $E$ & $n_{0}$ & $25$ & $36$ & $49$ & $64$ & $81$ \\ \hline
\multicolumn{1}{|l|}{$11a1$} & \multicolumn{1}{|l|}{$1$} &
\multicolumn{1}{|l|}{$0.312985$} &
\multicolumn{1}{|l|}{$0.179629$} &
\multicolumn{1}{|l|}{$0.239533$} &
\multicolumn{1}{|l|}{$0.144443$} &
\multicolumn{1}{|l|}{$0.235818$} \\ \hline
\multicolumn{1}{|l|}{$11a1$} & \multicolumn{1}{|l|}{$3$} &
\multicolumn{1}{|l|}{$0.423453$} &
\multicolumn{1}{|l|}{$0.208141$} & \multicolumn{1}{|l|}{$0.27803$}
& \multicolumn{1}{|l|}{$0.141689$} &
\multicolumn{1}{|l|}{$0.258489$} \\ \hline
\multicolumn{1}{|l|}{$11a1$} & \multicolumn{1}{|l|}{$5$} &
\multicolumn{1}{|l|}{$0.320314$} &
\multicolumn{1}{|l|}{$0.180445$} &
\multicolumn{1}{|l|}{$0.238248$} &
\multicolumn{1}{|l|}{$0.140963$} &
\multicolumn{1}{|l|}{$0.234221$} \\ \hline
\multicolumn{1}{|l|}{$11a1$} & \multicolumn{1}{|l|}{$15$} &
\multicolumn{1}{|l|}{$0.415875$} &
\multicolumn{1}{|l|}{$0.207029$} &
\multicolumn{1}{|l|}{$0.278424$} &
\multicolumn{1}{|l|}{$0.146673$} &
\multicolumn{1}{|l|}{$0.249818$} \\ \hline
\multicolumn{1}{|l|}{$11a1$} & \multicolumn{1}{|l|}{$23$} &
\multicolumn{1}{|l|}{$0.411152$} & \multicolumn{1}{|l|}{$0.20807$}
& \multicolumn{1}{|l|}{$0.283056$} &
\multicolumn{1}{|l|}{$0.149102$} &
\multicolumn{1}{|l|}{$0.250205$} \\ \hline
\multicolumn{1}{|l|}{$11a1$} & \multicolumn{1}{|l|}{$31$} &
\multicolumn{1}{|l|}{$0.412866$} &
\multicolumn{1}{|l|}{$0.205186$} &
\multicolumn{1}{|l|}{$0.277763$} &
\multicolumn{1}{|l|}{$0.150807$} &
\multicolumn{1}{|l|}{$0.254808$} \\ \hline
\multicolumn{1}{|l|}{$11a1$} & \multicolumn{1}{|l|}{$37$} &
\multicolumn{1}{|l|}{$0.314151$} &
\multicolumn{1}{|l|}{$0.171081$} &
\multicolumn{1}{|l|}{$0.235285$} &
\multicolumn{1}{|l|}{$0.143295$} &
\multicolumn{1}{|l|}{$0.230991$} \\ \hline
\end{tabular}

\begin{tabular}{|c|c|c|c|c|c|c|}
\hline $E$ & $n_{0}$ & $100$ & $121$ & $144$ & $169$ & $196$ \\
\hline \multicolumn{1}{|l|}{$11a1$} & \multicolumn{1}{|l|}{$1$} &
\multicolumn{1}{|l|}{$0.149895$} &
\multicolumn{1}{|l|}{$0.188637$} &
\multicolumn{1}{|l|}{$0.115865$} &
\multicolumn{1}{|l|}{$0.171031$} &
\multicolumn{1}{|l|}{$0.091912$} \\ \hline
\multicolumn{1}{|l|}{$11a1$} & \multicolumn{1}{|l|}{$3$} &
\multicolumn{1}{|l|}{$0.144478$} &
\multicolumn{1}{|l|}{$0.144478$} &
\multicolumn{1}{|l|}{$0.099708$} &
\multicolumn{1}{|l|}{$0.158332$} &
\multicolumn{1}{|l|}{$0.068375$} \\ \hline
\multicolumn{1}{|l|}{$11a1$} & \multicolumn{1}{|l|}{$5$} &
\multicolumn{1}{|l|}{$0.151277$} &
\multicolumn{1}{|l|}{$0.183852$} &
\multicolumn{1}{|l|}{$0.117823$} &
\multicolumn{1}{|l|}{$0.165406$} &
\multicolumn{1}{|l|}{$0.089266$} \\ \hline
\multicolumn{1}{|l|}{$11a1$} & \multicolumn{1}{|l|}{$15$} &
\multicolumn{1}{|l|}{$0.1391$} & \multicolumn{1}{|l|}{$0.1391$} &
\multicolumn{1}{|l|}{$0.096762$} &
\multicolumn{1}{|l|}{$0.156685$} &
\multicolumn{1}{|l|}{$0.069026$} \\ \hline
\multicolumn{1}{|l|}{$11a1$} & \multicolumn{1}{|l|}{$23$} &
\multicolumn{1}{|l|}{$0.140066$} &
\multicolumn{1}{|l|}{$0.190515$} & \multicolumn{1}{|l|}{$0.09152$}
& \multicolumn{1}{|l|}{$0.159746$} &
\multicolumn{1}{|l|}{$0.070316$} \\ \hline
\multicolumn{1}{|l|}{$11a1$} & \multicolumn{1}{|l|}{$31$} &
\multicolumn{1}{|l|}{$0.141375$} &
\multicolumn{1}{|l|}{$0.185663$} &
\multicolumn{1}{|l|}{$0.098614$} &
\multicolumn{1}{|l|}{$0.162711$} &
\multicolumn{1}{|l|}{$0.071718$} \\ \hline
\multicolumn{1}{|l|}{$11a1$} & \multicolumn{1}{|l|}{$37$} &
\multicolumn{1}{|l|}{$0.152059$} &
\multicolumn{1}{|l|}{$0.184835$} &
\multicolumn{1}{|l|}{$0.112431$} &
\multicolumn{1}{|l|}{$0.172438$} & \multicolumn{1}{|l|}{$0.08812$}
\\ \hline
\end{tabular}

\begin{tabular}{|c|c|c|c|c|c|c|}
\hline $E$ & $n_{0}$ & $225$ & $256$ & $289$ & $324$ & $361$ \\
\hline \multicolumn{1}{|l|}{$11a1$} & \multicolumn{1}{|l|}{$1$} &
\multicolumn{1}{|l|}{$0.203897$} &
\multicolumn{1}{|l|}{$0.076341$} &
\multicolumn{1}{|l|}{$0.135676$} & \multicolumn{1}{|l|}{$0.07271$}
& \multicolumn{1}{|l|}{$0.117004$} \\ \hline
\multicolumn{1}{|l|}{$11a1$} & \multicolumn{1}{|l|}{$3$} &
\multicolumn{1}{|l|}{$0.178166$} &
\multicolumn{1}{|l|}{$0.054526$} &
\multicolumn{1}{|l|}{$0.113502$} &
\multicolumn{1}{|l|}{$0.047742$} &
\multicolumn{1}{|l|}{$0.097803$} \\ \hline
\multicolumn{1}{|l|}{$11a1$} & \multicolumn{1}{|l|}{$5$} &
\multicolumn{1}{|l|}{$0.201578$} &
\multicolumn{1}{|l|}{$0.076127$} &
\multicolumn{1}{|l|}{$0.132841$} &
\multicolumn{1}{|l|}{$0.070794$} &
\multicolumn{1}{|l|}{$0.116838$} \\ \hline
\multicolumn{1}{|l|}{$11a1$} & \multicolumn{1}{|l|}{$15$} &
\multicolumn{1}{|l|}{$0.185691$} &
\multicolumn{1}{|l|}{$0.055885$} &
\multicolumn{1}{|l|}{$0.117211$} & \multicolumn{1}{|l|}{$0.0487$}
& \multicolumn{1}{|l|}{$0.094476$} \\ \hline
\multicolumn{1}{|l|}{$11a1$} & \multicolumn{1}{|l|}{$23$} &
\multicolumn{1}{|l|}{$0.179255$} &
\multicolumn{1}{|l|}{$0.058757$} &
\multicolumn{1}{|l|}{$0.112958$} &
\multicolumn{1}{|l|}{$0.048998$} & \multicolumn{1}{|l|}{$0.0973$}
\\ \hline \multicolumn{1}{|l|}{$11a1$} &
\multicolumn{1}{|l|}{$31$} & \multicolumn{1}{|l|}{$0.18378$} &
\multicolumn{1}{|l|}{$0.054368$} &
\multicolumn{1}{|l|}{$0.116274$} &
\multicolumn{1}{|l|}{$0.049585$} & \multicolumn{1}{|l|}{$0.0968$}
\\ \hline \multicolumn{1}{|l|}{$11a1$} &
\multicolumn{1}{|l|}{$37$} & \multicolumn{1}{|l|}{$0.208334$} &
\multicolumn{1}{|l|}{$0.079288$} &
\multicolumn{1}{|l|}{$0.132626$} & \multicolumn{1}{|l|}{$0.06951$}
& \multicolumn{1}{|l|}{$0.117505$} \\ \hline
\end{tabular}

\begin{tabular}{|c|c|c|c|c|c|c|}
\hline $E$ & $n_{0}$ & $0$ & $1$ & $4$ & $9$ & $16$ \\ \hline
\multicolumn{1}{|l|}{$14a1$} & \multicolumn{1}{|l|}{$1$} &
\multicolumn{1}{|l|}{$0.283019$} &
\multicolumn{1}{|l|}{$0.386791$} &
\multicolumn{1}{|l|}{$0.349053$} &
\multicolumn{1}{|l|}{$0.485425$} &
\multicolumn{1}{|l|}{$0.289702$} \\ \hline
\multicolumn{1}{|l|}{$14a1$} & \multicolumn{1}{|l|}{$15$} &
\multicolumn{1}{|l|}{$0.319039$} &
\multicolumn{1}{|l|}{$0.483879$} &
\multicolumn{1}{|l|}{$0.409463$} &
\multicolumn{1}{|l|}{$0.559197$} &
\multicolumn{1}{|l|}{$0.319645$} \\ \hline
\multicolumn{1}{|l|}{$14a1$} & \multicolumn{1}{|l|}{$23$} &
\multicolumn{1}{|l|}{$0.336754$} &
\multicolumn{1}{|l|}{$0.461198$} &
\multicolumn{1}{|l|}{$0.402525$} &
\multicolumn{1}{|l|}{$0.567646$} &
\multicolumn{1}{|l|}{$0.312323$} \\ \hline
\multicolumn{1}{|l|}{$14a1$} & \multicolumn{1}{|l|}{$29$} &
\multicolumn{1}{|l|}{$0.442312$} &
\multicolumn{1}{|l|}{$0.172938$} &
\multicolumn{1}{|l|}{$0.560877$} &
\multicolumn{1}{|l|}{$0.194746$} &
\multicolumn{1}{|l|}{$0.485192$} \\ \hline
\multicolumn{1}{|l|}{$14a1$} & \multicolumn{1}{|l|}{$37$} &
\multicolumn{1}{|l|}{$0.42059$} & \multicolumn{1}{|l|}{$0.175757$}
& \multicolumn{1}{|l|}{$0.589686$} &
\multicolumn{1}{|l|}{$0.185407$} &
\multicolumn{1}{|l|}{$0.492192$} \\ \hline
\multicolumn{1}{|l|}{$14a1$} & \multicolumn{1}{|l|}{$39$} &
\multicolumn{1}{|l|}{$0.312339$} &
\multicolumn{1}{|l|}{$0.493768$} &
\multicolumn{1}{|l|}{$0.407928$} &
\multicolumn{1}{|l|}{$0.536247$} &
\multicolumn{1}{|l|}{$0.328624$} \\ \hline
\multicolumn{1}{|l|}{$14a1$} & \multicolumn{1}{|l|}{$53$} &
\multicolumn{1}{|l|}{$0.447676$} &
\multicolumn{1}{|l|}{$0.171374$} &
\multicolumn{1}{|l|}{$0.571985$} &
\multicolumn{1}{|l|}{$0.191641$} &
\multicolumn{1}{|l|}{$0.472537$} \\ \hline
\end{tabular}

\begin{tabular}{|c|c|c|c|c|c|c|}
\hline $E$ & $n_{0}$ & $25$ & $36$ & $49$ & $64$ & $81$ \\ \hline
\multicolumn{1}{|l|}{$14a1$} & \multicolumn{1}{|l|}{$1$} &
\multicolumn{1}{|l|}{$0.28595$} & \multicolumn{1}{|l|}{$0.326485$}
& \multicolumn{1}{|l|}{$0.190388$} &
\multicolumn{1}{|l|}{$0.148085$} &
\multicolumn{1}{|l|}{$0.292196$} \\ \hline
\multicolumn{1}{|l|}{$14a1$} & \multicolumn{1}{|l|}{$15$} &
\multicolumn{1}{|l|}{$0.314341$} &
\multicolumn{1}{|l|}{$0.322687$} &
\multicolumn{1}{|l|}{$0.235968$} &
\multicolumn{1}{|l|}{$0.181454$} &
\multicolumn{1}{|l|}{$0.290788$} \\ \hline
\multicolumn{1}{|l|}{$14a1$} & \multicolumn{1}{|l|}{$23$} &
\multicolumn{1}{|l|}{$0.314975$} &
\multicolumn{1}{|l|}{$0.331965$} &
\multicolumn{1}{|l|}{$0.236908$} &
\multicolumn{1}{|l|}{$0.173742$} &
\multicolumn{1}{|l|}{$0.299374$} \\ \hline
\multicolumn{1}{|l|}{$14a1$} & \multicolumn{1}{|l|}{$29$} &
\multicolumn{1}{|l|}{$0.099951$} &
\multicolumn{1}{|l|}{$0.567019$} &
\multicolumn{1}{|l|}{$0.071551$} &
\multicolumn{1}{|l|}{$0.310347$} &
\multicolumn{1}{|l|}{$0.086056$} \\ \hline
\multicolumn{1}{|l|}{$14a1$} & \multicolumn{1}{|l|}{$37$} &
\multicolumn{1}{|l|}{$0.108121$} &
\multicolumn{1}{|l|}{$0.544933$} &
\multicolumn{1}{|l|}{$0.076132$} & \multicolumn{1}{|l|}{$0.32644$}
& \multicolumn{1}{|l|}{$0.081628$} \\ \hline
\multicolumn{1}{|l|}{$14a1$} & \multicolumn{1}{|l|}{$39$} &
\multicolumn{1}{|l|}{$0.327089$} &
\multicolumn{1}{|l|}{$0.316292$} &
\multicolumn{1}{|l|}{$0.251003$} &
\multicolumn{1}{|l|}{$0.186261$} & \multicolumn{1}{|l|}{$0.27054$}
\\ \hline \multicolumn{1}{|l|}{$14a1$} &
\multicolumn{1}{|l|}{$53$} & \multicolumn{1}{|l|}{$0.104698$} &
\multicolumn{1}{|l|}{$0.561543$} &
\multicolumn{1}{|l|}{$0.073051$} &
\multicolumn{1}{|l|}{$0.321251$} &
\multicolumn{1}{|l|}{$0.086696$} \\ \hline
\end{tabular}

\begin{tabular}{|c|c|c|c|c|c|c|}
\hline $E$ & $n_{0}$ & $100$ & $121$ & $144$ & $169$ & $196$ \\
\hline \multicolumn{1}{|l|}{$14a1$} & \multicolumn{1}{|l|}{$1$} &
\multicolumn{1}{|l|}{$0.142119$} &
\multicolumn{1}{|l|}{$0.149737$} &
\multicolumn{1}{|l|}{$0.159668$} &
\multicolumn{1}{|l|}{$0.119627$} &
\multicolumn{1}{|l|}{$0.087451$} \\ \hline
\multicolumn{1}{|l|}{$14a1$} & \multicolumn{1}{|l|}{$15$} &
\multicolumn{1}{|l|}{$0.139186$} &
\multicolumn{1}{|l|}{$0.150902$} &
\multicolumn{1}{|l|}{$0.143504$} &
\multicolumn{1}{|l|}{$0.116318$} &
\multicolumn{1}{|l|}{$0.076208$} \\ \hline
\multicolumn{1}{|l|}{$14a1$} & \multicolumn{1}{|l|}{$23$} &
\multicolumn{1}{|l|}{$0.134807$} &
\multicolumn{1}{|l|}{$0.144059$} &
\multicolumn{1}{|l|}{$0.150461$} &
\multicolumn{1}{|l|}{$0.112789$} &
\multicolumn{1}{|l|}{$0.072932$} \\ \hline
\multicolumn{1}{|l|}{$14a1$} & \multicolumn{1}{|l|}{$29$} &
\multicolumn{1}{|l|}{$0.265114$} &
\multicolumn{1}{|l|}{$0.039297$} &
\multicolumn{1}{|l|}{$0.303978$} & \multicolumn{1}{|l|}{$0.03107$}
& \multicolumn{1}{|l|}{$0.174705$} \\ \hline
\multicolumn{1}{|l|}{$14a1$} & \multicolumn{1}{|l|}{$37$} &
\multicolumn{1}{|l|}{$0.276584$} &
\multicolumn{1}{|l|}{$0.042321$} &
\multicolumn{1}{|l|}{$0.285788$} &
\multicolumn{1}{|l|}{$0.029985$} & \multicolumn{1}{|l|}{$0.18166$}
\\ \hline \multicolumn{1}{|l|}{$14a1$} &
\multicolumn{1}{|l|}{$39$} & \multicolumn{1}{|l|}{$0.140236$} &
\multicolumn{1}{|l|}{$0.14836$} & \multicolumn{1}{|l|}{$0.135962$}
& \multicolumn{1}{|l|}{$0.11914$} &
\multicolumn{1}{|l|}{$0.079435$} \\ \hline
\multicolumn{1}{|l|}{$14a1$} & \multicolumn{1}{|l|}{$53$} &
\multicolumn{1}{|l|}{$0.264155$} &
\multicolumn{1}{|l|}{$0.036761$} &
\multicolumn{1}{|l|}{$0.306407$} &
\multicolumn{1}{|l|}{$0.029676$} &
\multicolumn{1}{|l|}{$0.174854$} \\ \hline
\end{tabular}

\begin{tabular}{|c|c|c|c|c|c|c|}
\hline $E$ & $n_{0}$ & $225$ & $256$ & $289$ & $324$ & $361$ \\
\hline \multicolumn{1}{|l|}{$14a1$} & \multicolumn{1}{|l|}{$1$} &
\multicolumn{1}{|l|}{$0.147626$} &
\multicolumn{1}{|l|}{$0.073324$} &
\multicolumn{1}{|l|}{$0.081276$} &
\multicolumn{1}{|l|}{$0.090333$} &
\multicolumn{1}{|l|}{$0.069123$} \\ \hline
\multicolumn{1}{|l|}{$14a1$} & \multicolumn{1}{|l|}{$15$} &
\multicolumn{1}{|l|}{$0.13048$} & \multicolumn{1}{|l|}{$0.060237$}
& \multicolumn{1}{|l|}{$0.071448$} &
\multicolumn{1}{|l|}{$0.066382$} &
\multicolumn{1}{|l|}{$0.058287$} \\ \hline
\multicolumn{1}{|l|}{$14a1$} & \multicolumn{1}{|l|}{$23$} &
\multicolumn{1}{|l|}{$0.138795$} &
\multicolumn{1}{|l|}{$0.059647$} &
\multicolumn{1}{|l|}{$0.070181$} &
\multicolumn{1}{|l|}{$0.074594$} &
\multicolumn{1}{|l|}{$0.054617$} \\ \hline
\multicolumn{1}{|l|}{$14a1$} & \multicolumn{1}{|l|}{$29$} &
\multicolumn{1}{|l|}{$0.033083$} &
\multicolumn{1}{|l|}{$0.144352$} &
\multicolumn{1}{|l|}{$0.017318$} &
\multicolumn{1}{|l|}{$0.190917$} &
\multicolumn{1}{|l|}{$0.011085$} \\ \hline
\multicolumn{1}{|l|}{$14a1$} & \multicolumn{1}{|l|}{$37$} &
\multicolumn{1}{|l|}{$0.031929$} &
\multicolumn{1}{|l|}{$0.141464$} & \multicolumn{1}{|l|}{$0.01651$}
& \multicolumn{1}{|l|}{$0.180291$} &
\multicolumn{1}{|l|}{$0.011341$} \\ \hline
\multicolumn{1}{|l|}{$14a1$} & \multicolumn{1}{|l|}{$39$} &
\multicolumn{1}{|l|}{$0.129477$} &
\multicolumn{1}{|l|}{$0.062409$} &
\multicolumn{1}{|l|}{$0.074189$} & \multicolumn{1}{|l|}{$0.06682$}
& \multicolumn{1}{|l|}{$0.059717$} \\ \hline
\multicolumn{1}{|l|}{$14a1$} & \multicolumn{1}{|l|}{$53$} &
\multicolumn{1}{|l|}{$0.03286$} & \multicolumn{1}{|l|}{$0.137962$}
& \multicolumn{1}{|l|}{$0.01691$} &
\multicolumn{1}{|l|}{$0.194275$} &
\multicolumn{1}{|l|}{$0.012076$} \\ \hline
\end{tabular}

\begin{tabular}{|l|l|l|l|l|l|l|}
\hline $E$ & $n_{0}$ & $0$ & $1$ & $4$ & $9$ & $16$ \\ \hline
$17a1$ & $3$ & $0.337432$ & $0.477285$ & $0.501361$ & $0.41195$ &
$0.402614$
\\ \hline
$17a1$ & $7$ & $0.333173$ & $0.480345$ & $0.512958$ & $0.411449$ &
$0.397752$
\\ \hline
$17a1$ & $11$ & $0.331548$ & $0.470597$ & $0.506991$ & $0.41595$ &
$0.398308$
\\ \hline
$17a1$ & $23$ & $0.324727$ & $0.469987$ & $0.510093$ & $0.413703$ & $%
0.409981 $ \\ \hline
$17a1$ & $31$ & $0.332091$ & $0.482396$ & $0.496191$ & $0.410295$ & $%
0.405293 $ \\ \hline $17a1$ & $39$ & $0.335686$ & $0.481485$ &
$0.50061$ & $0.403566$ & $0.403538$
\\ \hline
\end{tabular}

\begin{tabular}{|l|l|l|l|l|l|l|}
\hline $E$ & $n_{0}$ & $25$ & $36$ & $49$ & $64$ & $81$ \\ \hline
$17a1$ & $3$ & $0.291697$ & $0.298922$ & $0.227993$ & $0.217491$ &
$0.190683$
\\ \hline
$17a1$ & $7$ & $0.294852$ & $0.305199$ & $0.224537$ & $0.209766$ &
$0.191861$
\\ \hline
$17a1$ & $11$ & $0.29197$ & $0.307093$ & $0.224453$ & $0.212266$ &
$0.194131$
\\ \hline
$17a1$ & $23$ & $0.302838$ & $0.296815$ & $0.223042$ & $0.212759$ & $%
0.197424 $ \\ \hline $17a1$ & $31$ & $0.299154$ & $0.303459$ &
$0.219757$ & $0.213895$ & $0.19307$
\\ \hline
$17a1$ & $39$ & $0.29654$ & $0.302868$ & $0.21993$ & $0.218869$ &
$0.19364$
\\ \hline
\end{tabular}

\begin{tabular}{|c|c|c|c|c|c|c|}
\hline $E$ & $n_{0}$ & $100$ & $121$ & $144$ & $169$ & $196$ \\
\hline \multicolumn{1}{|l|}{$17a1$} & \multicolumn{1}{|l|}{$3$} &
\multicolumn{1}{|l|}{$0.155873$} &
\multicolumn{1}{|l|}{$0.138493$} &
\multicolumn{1}{|l|}{$0.129301$} &
\multicolumn{1}{|l|}{$0.109835$} &
\multicolumn{1}{|l|}{$0.086544$} \\ \hline
\multicolumn{1}{|l|}{$17a1$} & \multicolumn{1}{|l|}{$7$} &
\multicolumn{1}{|l|}{$0.153253$} &
\multicolumn{1}{|l|}{$0.134025$} &
\multicolumn{1}{|l|}{$0.127087$} &
\multicolumn{1}{|l|}{$0.109443$} & \multicolumn{1}{|l|}{$0.08902$}
\\ \hline \multicolumn{1}{|l|}{$17a1$} &
\multicolumn{1}{|l|}{$11$} & \multicolumn{1}{|l|}{$0.152146$} &
\multicolumn{1}{|l|}{$0.142299$} & \multicolumn{1}{|l|}{$0.1244$}
& \multicolumn{1}{|l|}{$0.115795$} &
\multicolumn{1}{|l|}{$0.090622$} \\ \hline
\multicolumn{1}{|l|}{$17a1$} & \multicolumn{1}{|l|}{$23$} &
\multicolumn{1}{|l|}{$0.149471$} &
\multicolumn{1}{|l|}{$0.136131$} &
\multicolumn{1}{|l|}{$0.125797$} &
\multicolumn{1}{|l|}{$0.115204$} &
\multicolumn{1}{|l|}{$0.088695$} \\ \hline
\multicolumn{1}{|l|}{$17a1$} & \multicolumn{1}{|l|}{$31$} &
\multicolumn{1}{|l|}{$0.153817$} &
\multicolumn{1}{|l|}{$0.141537$} &
\multicolumn{1}{|l|}{$0.128656$} & \multicolumn{1}{|l|}{$0.10723$}
& \multicolumn{1}{|l|}{$0.086032$} \\ \hline
\multicolumn{1}{|l|}{$17a1$} & \multicolumn{1}{|l|}{$39$} &
\multicolumn{1}{|l|}{$0.155061$} &
\multicolumn{1}{|l|}{$0.139321$} & \multicolumn{1}{|l|}{$0.12622$}
& \multicolumn{1}{|l|}{$0.110953$} &
\multicolumn{1}{|l|}{$0.084326$} \\ \hline
\end{tabular}

\begin{tabular}{|c|c|c|c|c|c|c|}
\hline $E$ & $n_{0}$ & $225$ & $256$ & $289$ & $324$ & $361$ \\
\hline \multicolumn{1}{|l|}{$17a1$} & \multicolumn{1}{|l|}{$3$} &
\multicolumn{1}{|l|}{$0.10133$} & \multicolumn{1}{|l|}{$0.066479$}
& \multicolumn{1}{|l|}{$0.070144$} &
\multicolumn{1}{|l|}{$0.054392$} &
\multicolumn{1}{|l|}{$0.063161$} \\ \hline
\multicolumn{1}{|l|}{$17a1$} & \multicolumn{1}{|l|}{$7$} &
\multicolumn{1}{|l|}{$0.102284$} &
\multicolumn{1}{|l|}{$0.066453$} & \multicolumn{1}{|l|}{$0.07182$}
& \multicolumn{1}{|l|}{$0.052183$} &
\multicolumn{1}{|l|}{$0.061862$} \\ \hline
\multicolumn{1}{|l|}{$17a1$} & \multicolumn{1}{|l|}{$11$} &
\multicolumn{1}{|l|}{$0.101815$} & \multicolumn{1}{|l|}{$0.06621$}
& \multicolumn{1}{|l|}{$0.073766$} &
\multicolumn{1}{|l|}{$0.053459$} &
\multicolumn{1}{|l|}{$0.060803$} \\ \hline
\multicolumn{1}{|l|}{$17a1$} & \multicolumn{1}{|l|}{$23$} &
\multicolumn{1}{|l|}{$0.104502$} &
\multicolumn{1}{|l|}{$0.066214$} &
\multicolumn{1}{|l|}{$0.069106$} &
\multicolumn{1}{|l|}{$0.055098$} &
\multicolumn{1}{|l|}{$0.060021$} \\ \hline
\multicolumn{1}{|l|}{$17a1$} & \multicolumn{1}{|l|}{$31$} &
\multicolumn{1}{|l|}{$0.106572$} &
\multicolumn{1}{|l|}{$0.068254$} &
\multicolumn{1}{|l|}{$0.071002$} &
\multicolumn{1}{|l|}{$0.057448$} &
\multicolumn{1}{|l|}{$0.058929$} \\ \hline
\multicolumn{1}{|l|}{$17a1$} & \multicolumn{1}{|l|}{$39$} &
\multicolumn{1}{|l|}{$0.108346$} &
\multicolumn{1}{|l|}{$0.065278$} &
\multicolumn{1}{|l|}{$0.073044$} &
\multicolumn{1}{|l|}{$0.055616$} &
\multicolumn{1}{|l|}{$0.058537$} \\ \hline
\end{tabular}

\begin{tabular}{|c|c|c|c|c|c|c|}
\hline $E$ & $n_{0}$ & $0$ & $1$ & $4$ & $9$ & $16$ \\ \hline
\multicolumn{1}{|l|}{$20a1$} & \multicolumn{1}{|l|}{$1$} &
\multicolumn{1}{|l|}{$0.268253$} & \multicolumn{1}{|l|}{$0.3465$}
& \multicolumn{1}{|l|}{$0.315475$} &
\multicolumn{1}{|l|}{$0.427111$} & \multicolumn{1}{|l|}{$0.27129$}
\\ \hline \multicolumn{1}{|l|}{$20a1$} &
\multicolumn{1}{|l|}{$21$} & \multicolumn{1}{|l|}{$0.266462$} &
\multicolumn{1}{|l|}{$0.337876$} & \multicolumn{1}{|l|}{$0.32056$}
& \multicolumn{1}{|l|}{$0.431508$} &
\multicolumn{1}{|l|}{$0.272359$} \\ \hline
\multicolumn{1}{|l|}{$20a1$} & \multicolumn{1}{|l|}{$29$} &
\multicolumn{1}{|l|}{$0.267792$} &
\multicolumn{1}{|l|}{$0.343135$} &
\multicolumn{1}{|l|}{$0.317666$} &
\multicolumn{1}{|l|}{$0.425236$} &
\multicolumn{1}{|l|}{$0.271235$} \\ \hline
\end{tabular}

\begin{tabular}{|c|c|c|c|c|c|c|}
\hline $E$ & $n_{0}$ & $25$ & $36$ & $49$ & $64$ & $81$ \\ \hline
\multicolumn{1}{|l|}{$20a1$} & \multicolumn{1}{|l|}{$1$} &
\multicolumn{1}{|l|}{$0.254326$} &
\multicolumn{1}{|l|}{$0.307296$} &
\multicolumn{1}{|l|}{$0.210761$} &
\multicolumn{1}{|l|}{$0.179752$} &
\multicolumn{1}{|l|}{$0.245513$} \\ \hline
\multicolumn{1}{|l|}{$20a1$} & \multicolumn{1}{|l|}{$21$} &
\multicolumn{1}{|l|}{$0.253463$} &
\multicolumn{1}{|l|}{$0.304903$} &
\multicolumn{1}{|l|}{$0.209301$} &
\multicolumn{1}{|l|}{$0.178783$} &
\multicolumn{1}{|l|}{$0.246748$} \\ \hline
\multicolumn{1}{|l|}{$20a1$} & \multicolumn{1}{|l|}{$29$} &
\multicolumn{1}{|l|}{$0.258567$} &
\multicolumn{1}{|l|}{$0.308143$} & \multicolumn{1}{|l|}{$0.20873$}
& \multicolumn{1}{|l|}{$0.178674$} &
\multicolumn{1}{|l|}{$0.252364$} \\ \hline
\end{tabular}

\begin{tabular}{|c|c|c|c|c|c|c|}
\hline $E$ & $n_{0}$ & $100$ & $121$ & $144$ & $169$ & $196$ \\
\hline \multicolumn{1}{|l|}{$20a1$} & \multicolumn{1}{|l|}{$1$} &
\multicolumn{1}{|l|}{$0.144222$} &
\multicolumn{1}{|l|}{$0.141449$} &
\multicolumn{1}{|l|}{$0.171656$} &
\multicolumn{1}{|l|}{$0.115768$} &
\multicolumn{1}{|l|}{$0.095252$} \\ \hline
\multicolumn{1}{|l|}{$20a1$} & \multicolumn{1}{|l|}{$21$} &
\multicolumn{1}{|l|}{$0.146098$} &
\multicolumn{1}{|l|}{$0.144796$} &
\multicolumn{1}{|l|}{$0.165373$} &
\multicolumn{1}{|l|}{$0.115249$} & \multicolumn{1}{|l|}{$0.09443$}
\\ \hline \multicolumn{1}{|l|}{$20a1$} &
\multicolumn{1}{|l|}{$29$} & \multicolumn{1}{|l|}{$0.142937$} &
\multicolumn{1}{|l|}{$0.14178$} & \multicolumn{1}{|l|}{$0.1634$} &
\multicolumn{1}{|l|}{$0.115021$} & \multicolumn{1}{|l|}{$0.09569$}
\\ \hline
\end{tabular}

\begin{tabular}{|c|c|c|c|c|c|c|}
\hline $E$ & $n_{0}$ & $225$ & $256$ & $289$ & $324$ & $361$ \\
\hline \multicolumn{1}{|l|}{$20a1$} & \multicolumn{1}{|l|}{$1$} &
\multicolumn{1}{|l|}{$0.141739$} &
\multicolumn{1}{|l|}{$0.076533$} &
\multicolumn{1}{|l|}{$0.082182$} &
\multicolumn{1}{|l|}{$0.091834$} &
\multicolumn{1}{|l|}{$0.066588$} \\ \hline
\multicolumn{1}{|l|}{$20a1$} & \multicolumn{1}{|l|}{$21$} &
\multicolumn{1}{|l|}{$0.147373$} &
\multicolumn{1}{|l|}{$0.078015$} &
\multicolumn{1}{|l|}{$0.082924$} &
\multicolumn{1}{|l|}{$0.091549$} &
\multicolumn{1}{|l|}{$0.067251$} \\ \hline
\multicolumn{1}{|l|}{$20a1$} & \multicolumn{1}{|l|}{$29$} &
\multicolumn{1}{|l|}{$0.14228$} & \multicolumn{1}{|l|}{$0.081021$}
& \multicolumn{1}{|l|}{$0.081558$} &
\multicolumn{1}{|l|}{$0.091311$} &
\multicolumn{1}{|l|}{$0.067867$} \\ \hline
\end{tabular}

\begin{tabular}{|c|c|c|c|c|c|c|}
\hline $E$ & $n_{0}$ & $0$ & $1$ & $4$ & $9$ & $16$ \\ \hline
\multicolumn{1}{|l|}{$34a1$} & \multicolumn{1}{|l|}{$1$} &
\multicolumn{1}{|l|}{$0.300968$} &
\multicolumn{1}{|l|}{$0.385865$} &
\multicolumn{1}{|l|}{$0.387258$} &
\multicolumn{1}{|l|}{$0.462396$} & \multicolumn{1}{|l|}{$0.28402$}
\\ \hline \multicolumn{1}{|l|}{$34a1$} &
\multicolumn{1}{|l|}{$13$} & \multicolumn{1}{|l|}{$0.290206$} &
\multicolumn{1}{|l|}{$0.415303$} &
\multicolumn{1}{|l|}{$0.352209$} &
\multicolumn{1}{|l|}{$0.474225$} &
\multicolumn{1}{|l|}{$0.272241$} \\ \hline
\multicolumn{1}{|l|}{$34a1$} & \multicolumn{1}{|l|}{$19$} &
\multicolumn{1}{|l|}{$0.353435$} &
\multicolumn{1}{|l|}{$0.475157$} &
\multicolumn{1}{|l|}{$0.436218$} &
\multicolumn{1}{|l|}{$0.505167$} &
\multicolumn{1}{|l|}{$0.317592$} \\ \hline
\multicolumn{1}{|l|}{$34a1$} & \multicolumn{1}{|l|}{$21$} &
\multicolumn{1}{|l|}{$0.29167$} & \multicolumn{1}{|l|}{$0.415613$}
& \multicolumn{1}{|l|}{$0.359182$} &
\multicolumn{1}{|l|}{$0.472539$} &
\multicolumn{1}{|l|}{$0.274045$} \\ \hline
\multicolumn{1}{|l|}{$34a1$} & \multicolumn{1}{|l|}{$33$} &
\multicolumn{1}{|l|}{$0.30458$} & \multicolumn{1}{|l|}{$0.388798$}
& \multicolumn{1}{|l|}{$0.381037$} &
\multicolumn{1}{|l|}{$0.440285$} &
\multicolumn{1}{|l|}{$0.291886$} \\ \hline
\multicolumn{1}{|l|}{$34a1$} & \multicolumn{1}{|l|}{$35$} &
\multicolumn{1}{|l|}{$0.357437$} & \multicolumn{1}{|l|}{$0.47347$}
& \multicolumn{1}{|l|}{$0.42035$} &
\multicolumn{1}{|l|}{$0.504132$} &
\multicolumn{1}{|l|}{$0.326558$} \\ \hline
\multicolumn{1}{|l|}{$34a1$} & \multicolumn{1}{|l|}{$43$} &
\multicolumn{1}{|l|}{$0.355486$} &
\multicolumn{1}{|l|}{$0.466179$} & \multicolumn{1}{|l|}{$0.44077$}
& \multicolumn{1}{|l|}{$0.503479$} &
\multicolumn{1}{|l|}{$0.323861$} \\ \hline
\multicolumn{1}{|l|}{$34a1$} & \multicolumn{1}{|l|}{$53$} &
\multicolumn{1}{|l|}{$0.281215$} &
\multicolumn{1}{|l|}{$0.413834$} &
\multicolumn{1}{|l|}{$0.357105$} &
\multicolumn{1}{|l|}{$0.470171$} &
\multicolumn{1}{|l|}{$0.283536$} \\ \hline
\multicolumn{1}{|l|}{$34a1$} & \multicolumn{1}{|l|}{$59$} &
\multicolumn{1}{|l|}{$0.345971$} &
\multicolumn{1}{|l|}{$0.471297$} &
\multicolumn{1}{|l|}{$0.436144$} & \multicolumn{1}{|l|}{$0.50406$}
& \multicolumn{1}{|l|}{$0.327326$} \\ \hline
\multicolumn{1}{|l|}{$34a1$} & \multicolumn{1}{|l|}{$67$} &
\multicolumn{1}{|l|}{$0.351665$} &
\multicolumn{1}{|l|}{$0.467335$} &
\multicolumn{1}{|l|}{$0.427308$} &
\multicolumn{1}{|l|}{$0.512714$} &
\multicolumn{1}{|l|}{$0.326024$} \\ \hline
\multicolumn{1}{|l|}{$34a1$} & \multicolumn{1}{|l|}{$69$} &
\multicolumn{1}{|l|}{$0.290429$} &
\multicolumn{1}{|l|}{$0.408492$} &
\multicolumn{1}{|l|}{$0.366839$} &
\multicolumn{1}{|l|}{$0.478386$} &
\multicolumn{1}{|l|}{$0.275193$} \\ \hline
\multicolumn{1}{|l|}{$34a1$} & \multicolumn{1}{|l|}{$77$} &
\multicolumn{1}{|l|}{$0.293768$} & \multicolumn{1}{|l|}{$0.41554$}
& \multicolumn{1}{|l|}{$0.350608$} &
\multicolumn{1}{|l|}{$0.478178$} &
\multicolumn{1}{|l|}{$0.272873$} \\ \hline
\multicolumn{1}{|l|}{$34a1$} & \multicolumn{1}{|l|}{$83$} &
\multicolumn{1}{|l|}{$0.352644$} &
\multicolumn{1}{|l|}{$0.475251$} &
\multicolumn{1}{|l|}{$0.440215$} &
\multicolumn{1}{|l|}{$0.500611$} &
\multicolumn{1}{|l|}{$0.328119$} \\ \hline
\multicolumn{1}{|l|}{$34a1$} & \multicolumn{1}{|l|}{$89$} &
\multicolumn{1}{|l|}{$0.305732$} &
\multicolumn{1}{|l|}{$0.396955$} &
\multicolumn{1}{|l|}{$0.372179$} &
\multicolumn{1}{|l|}{$0.443078$} &
\multicolumn{1}{|l|}{$0.296228$} \\ \hline
\multicolumn{1}{|l|}{$34a1$} & \multicolumn{1}{|l|}{$93$} &
\multicolumn{1}{|l|}{$0.283804$} & \multicolumn{1}{|l|}{$0.42395$}
& \multicolumn{1}{|l|}{$0.358696$} &
\multicolumn{1}{|l|}{$0.479956$} &
\multicolumn{1}{|l|}{$0.279951$} \\ \hline
\multicolumn{1}{|l|}{$34a1$} & \multicolumn{1}{|l|}{$101$} &
\multicolumn{1}{|l|}{$0.29705$} & \multicolumn{1}{|l|}{$0.412981$}
& \multicolumn{1}{|l|}{$0.359811$} &
\multicolumn{1}{|l|}{$0.476887$} &
\multicolumn{1}{|l|}{$0.286045$} \\ \hline
\multicolumn{1}{|l|}{$34a1$} & \multicolumn{1}{|l|}{$115$} &
\multicolumn{1}{|l|}{$0.34747$} & \multicolumn{1}{|l|}{$0.476572$}
& \multicolumn{1}{|l|}{$0.438912$} &
\multicolumn{1}{|l|}{$0.505538$} &
\multicolumn{1}{|l|}{$0.321909$} \\ \hline
\multicolumn{1}{|l|}{$34a1$} & \multicolumn{1}{|l|}{$117$} &
\multicolumn{1}{|l|}{$0.291683$} &
\multicolumn{1}{|l|}{$0.420945$} &
\multicolumn{1}{|l|}{$0.355004$} &
\multicolumn{1}{|l|}{$0.476725$} &
\multicolumn{1}{|l|}{$0.278145$} \\ \hline
\multicolumn{1}{|l|}{$34a1$} & \multicolumn{1}{|l|}{$123$} &
\multicolumn{1}{|l|}{$0.354215$} &
\multicolumn{1}{|l|}{$0.475638$} &
\multicolumn{1}{|l|}{$0.437478$} &
\multicolumn{1}{|l|}{$0.495364$} & \multicolumn{1}{|l|}{$0.32921$}
\\ \hline
\end{tabular}

\begin{tabular}{|c|c|c|c|c|c|c|}
\hline $E$ & $n_{0}$ & $25$ & $36$ & $49$ & $64$ & $81$ \\ \hline
\multicolumn{1}{|l|}{$34a1$} & \multicolumn{1}{|l|}{$1$} &
\multicolumn{1}{|l|}{$0.247423$} &
\multicolumn{1}{|l|}{$0.309932$} &
\multicolumn{1}{|l|}{$0.194351$} &
\multicolumn{1}{|l|}{$0.177262$} &
\multicolumn{1}{|l|}{$0.225383$} \\ \hline
\multicolumn{1}{|l|}{$34a1$} & \multicolumn{1}{|l|}{$13$} &
\multicolumn{1}{|l|}{$0.271392$} &
\multicolumn{1}{|l|}{$0.294956$} &
\multicolumn{1}{|l|}{$0.199301$} &
\multicolumn{1}{|l|}{$0.166857$} &
\multicolumn{1}{|l|}{$0.230411$} \\ \hline
\multicolumn{1}{|l|}{$34a1$} & \multicolumn{1}{|l|}{$19$} &
\multicolumn{1}{|l|}{$0.271006$} &
\multicolumn{1}{|l|}{$0.324198$} &
\multicolumn{1}{|l|}{$0.191454$} &
\multicolumn{1}{|l|}{$0.171407$} &
\multicolumn{1}{|l|}{$0.214279$} \\ \hline
\multicolumn{1}{|l|}{$34a1$} & \multicolumn{1}{|l|}{$21$} &
\multicolumn{1}{|l|}{$0.265973$} &
\multicolumn{1}{|l|}{$0.301452$} & \multicolumn{1}{|l|}{$0.19956$}
& \multicolumn{1}{|l|}{$0.159942$} &
\multicolumn{1}{|l|}{$0.230879$} \\ \hline
\multicolumn{1}{|l|}{$34a1$} & \multicolumn{1}{|l|}{$33$} &
\multicolumn{1}{|l|}{$0.267835$} &
\multicolumn{1}{|l|}{$0.311831$} &
\multicolumn{1}{|l|}{$0.193117$} &
\multicolumn{1}{|l|}{$0.172183$} &
\multicolumn{1}{|l|}{$0.229097$} \\ \hline
\multicolumn{1}{|l|}{$34a1$} & \multicolumn{1}{|l|}{$35$} &
\multicolumn{1}{|l|}{$0.275831$} &
\multicolumn{1}{|l|}{$0.327544$} &
\multicolumn{1}{|l|}{$0.189914$} & \multicolumn{1}{|l|}{$0.17779$}
& \multicolumn{1}{|l|}{$0.220039$} \\ \hline
\multicolumn{1}{|l|}{$34a1$} & \multicolumn{1}{|l|}{$43$} &
\multicolumn{1}{|l|}{$0.272699$} &
\multicolumn{1}{|l|}{$0.317971$} &
\multicolumn{1}{|l|}{$0.202658$} &
\multicolumn{1}{|l|}{$0.174474$} &
\multicolumn{1}{|l|}{$0.211269$} \\ \hline
\multicolumn{1}{|l|}{$34a1$} & \multicolumn{1}{|l|}{$53$} &
\multicolumn{1}{|l|}{$0.277814$} &
\multicolumn{1}{|l|}{$0.310885$} &
\multicolumn{1}{|l|}{$0.201981$} &
\multicolumn{1}{|l|}{$0.161572$} &
\multicolumn{1}{|l|}{$0.229304$} \\ \hline
\multicolumn{1}{|l|}{$34a1$} & \multicolumn{1}{|l|}{$59$} &
\multicolumn{1}{|l|}{$0.267928$} &
\multicolumn{1}{|l|}{$0.327115$} &
\multicolumn{1}{|l|}{$0.193429$} &
\multicolumn{1}{|l|}{$0.175461$} &
\multicolumn{1}{|l|}{$0.215779$} \\ \hline
\multicolumn{1}{|l|}{$34a1$} & \multicolumn{1}{|l|}{$67$} &
\multicolumn{1}{|l|}{$0.273065$} &
\multicolumn{1}{|l|}{$0.324959$} &
\multicolumn{1}{|l|}{$0.192272$} &
\multicolumn{1}{|l|}{$0.172805$} &
\multicolumn{1}{|l|}{$0.212269$} \\ \hline
\multicolumn{1}{|l|}{$34a1$} & \multicolumn{1}{|l|}{$69$} &
\multicolumn{1}{|l|}{$0.267456$} & \multicolumn{1}{|l|}{$0.29777$}
& \multicolumn{1}{|l|}{$0.207129$} &
\multicolumn{1}{|l|}{$0.162831$} &
\multicolumn{1}{|l|}{$0.232521$} \\ \hline
\multicolumn{1}{|l|}{$34a1$} & \multicolumn{1}{|l|}{$77$} &
\multicolumn{1}{|l|}{$0.272458$} &
\multicolumn{1}{|l|}{$0.297814$} &
\multicolumn{1}{|l|}{$0.201069$} &
\multicolumn{1}{|l|}{$0.167444$} &
\multicolumn{1}{|l|}{$0.227964$} \\ \hline
\multicolumn{1}{|l|}{$34a1$} & \multicolumn{1}{|l|}{$83$} &
\multicolumn{1}{|l|}{$0.275118$} &
\multicolumn{1}{|l|}{$0.314132$} &
\multicolumn{1}{|l|}{$0.199511$} &
\multicolumn{1}{|l|}{$0.174251$} &
\multicolumn{1}{|l|}{$0.210163$} \\ \hline
\multicolumn{1}{|l|}{$34a1$} & \multicolumn{1}{|l|}{$89$} &
\multicolumn{1}{|l|}{$0.255453$} &
\multicolumn{1}{|l|}{$0.318219$} &
\multicolumn{1}{|l|}{$0.207944$} &
\multicolumn{1}{|l|}{$0.165966$} &
\multicolumn{1}{|l|}{$0.226874$} \\ \hline
\multicolumn{1}{|l|}{$34a1$} & \multicolumn{1}{|l|}{$93$} &
\multicolumn{1}{|l|}{$0.259963$} &
\multicolumn{1}{|l|}{$0.297419$} &
\multicolumn{1}{|l|}{$0.204218$} &
\multicolumn{1}{|l|}{$0.165222$} &
\multicolumn{1}{|l|}{$0.234308$} \\ \hline
\multicolumn{1}{|l|}{$34a1$} & \multicolumn{1}{|l|}{$101$} &
\multicolumn{1}{|l|}{$0.254251$} &
\multicolumn{1}{|l|}{$0.303388$} &
\multicolumn{1}{|l|}{$0.197465$} & \multicolumn{1}{|l|}{$0.15854$}
& \multicolumn{1}{|l|}{$0.23436$} \\ \hline
\multicolumn{1}{|l|}{$34a1$} & \multicolumn{1}{|l|}{$115$} &
\multicolumn{1}{|l|}{$0.270414$} &
\multicolumn{1}{|l|}{$0.323107$} &
\multicolumn{1}{|l|}{$0.196365$} &
\multicolumn{1}{|l|}{$0.177784$} & \multicolumn{1}{|l|}{$0.19878$}
\\ \hline \multicolumn{1}{|l|}{$34a1$} &
\multicolumn{1}{|l|}{$117$} & \multicolumn{1}{|l|}{$0.260653$} &
\multicolumn{1}{|l|}{$0.2887$} & \multicolumn{1}{|l|}{$0.202377$}
& \multicolumn{1}{|l|}{$0.157916$} &
\multicolumn{1}{|l|}{$0.244395$} \\ \hline
\multicolumn{1}{|l|}{$34a1$} & \multicolumn{1}{|l|}{$123$} &
\multicolumn{1}{|l|}{$0.270226$} &
\multicolumn{1}{|l|}{$0.335145$} &
\multicolumn{1}{|l|}{$0.200594$} &
\multicolumn{1}{|l|}{$0.170866$} &
\multicolumn{1}{|l|}{$0.208101$} \\ \hline
\end{tabular}

\begin{tabular}{|c|c|c|c|c|c|c|}
\hline $E$ & $n_{0}$ & $100$ & $121$ & $144$ & $169$ & $196$ \\
\hline \multicolumn{1}{|l|}{$34a1$} & \multicolumn{1}{|l|}{$1$} &
\multicolumn{1}{|l|}{$0.128304$} &
\multicolumn{1}{|l|}{$0.122231$} &
\multicolumn{1}{|l|}{$0.147564$} &
\multicolumn{1}{|l|}{$0.102664$} &
\multicolumn{1}{|l|}{$0.082786$} \\ \hline
\multicolumn{1}{|l|}{$34a1$} & \multicolumn{1}{|l|}{$13$} &
\multicolumn{1}{|l|}{$0.129592$} &
\multicolumn{1}{|l|}{$0.130917$} &
\multicolumn{1}{|l|}{$0.143784$} &
\multicolumn{1}{|l|}{$0.100675$} & \multicolumn{1}{|l|}{$0.07895$}
\\ \hline \multicolumn{1}{|l|}{$34a1$} &
\multicolumn{1}{|l|}{$19$} & \multicolumn{1}{|l|}{$0.130772$} &
\multicolumn{1}{|l|}{$0.109944$} &
\multicolumn{1}{|l|}{$0.140669$} &
\multicolumn{1}{|l|}{$0.086182$} &
\multicolumn{1}{|l|}{$0.077919$} \\ \hline
\multicolumn{1}{|l|}{$34a1$} & \multicolumn{1}{|l|}{$21$} &
\multicolumn{1}{|l|}{$0.132753$} &
\multicolumn{1}{|l|}{$0.126583$} &
\multicolumn{1}{|l|}{$0.143338$} &
\multicolumn{1}{|l|}{$0.106306$} &
\multicolumn{1}{|l|}{$0.080138$} \\ \hline
\multicolumn{1}{|l|}{$34a1$} & \multicolumn{1}{|l|}{$33$} &
\multicolumn{1}{|l|}{$0.131932$} &
\multicolumn{1}{|l|}{$0.123085$} &
\multicolumn{1}{|l|}{$0.140728$} &
\multicolumn{1}{|l|}{$0.099784$} &
\multicolumn{1}{|l|}{$0.082223$} \\ \hline
\multicolumn{1}{|l|}{$34a1$} & \multicolumn{1}{|l|}{$35$} &
\multicolumn{1}{|l|}{$0.130884$} &
\multicolumn{1}{|l|}{$0.111755$} &
\multicolumn{1}{|l|}{$0.141213$} & \multicolumn{1}{|l|}{$0.08385$}
& \multicolumn{1}{|l|}{$0.079269$} \\ \hline
\multicolumn{1}{|l|}{$34a1$} & \multicolumn{1}{|l|}{$43$} &
\multicolumn{1}{|l|}{$0.128261$} &
\multicolumn{1}{|l|}{$0.109375$} & \multicolumn{1}{|l|}{$0.13594$}
& \multicolumn{1}{|l|}{$0.083931$} &
\multicolumn{1}{|l|}{$0.071109$} \\ \hline
\multicolumn{1}{|l|}{$34a1$} & \multicolumn{1}{|l|}{$53$} &
\multicolumn{1}{|l|}{$0.12997$} & \multicolumn{1}{|l|}{$0.126922$}
& \multicolumn{1}{|l|}{$0.142644$} &
\multicolumn{1}{|l|}{$0.100092$} & \multicolumn{1}{|l|}{$0.07645$}
\\ \hline \multicolumn{1}{|l|}{$34a1$} &
\multicolumn{1}{|l|}{$59$} & \multicolumn{1}{|l|}{$0.137624$} &
\multicolumn{1}{|l|}{$0.106901$} &
\multicolumn{1}{|l|}{$0.140647$} & \multicolumn{1}{|l|}{$0.08392$}
& \multicolumn{1}{|l|}{$0.074337$} \\ \hline
\multicolumn{1}{|l|}{$34a1$} & \multicolumn{1}{|l|}{$67$} &
\multicolumn{1}{|l|}{$0.138161$} &
\multicolumn{1}{|l|}{$0.108241$} &
\multicolumn{1}{|l|}{$0.136687$} &
\multicolumn{1}{|l|}{$0.090009$} &
\multicolumn{1}{|l|}{$0.081168$} \\ \hline
\multicolumn{1}{|l|}{$34a1$} & \multicolumn{1}{|l|}{$69$} &
\multicolumn{1}{|l|}{$0.127983$} & \multicolumn{1}{|l|}{$0.12296$}
& \multicolumn{1}{|l|}{$0.146421$} &
\multicolumn{1}{|l|}{$0.100218$} &
\multicolumn{1}{|l|}{$0.072738$} \\ \hline
\multicolumn{1}{|l|}{$34a1$} & \multicolumn{1}{|l|}{$77$} &
\multicolumn{1}{|l|}{$0.130561$} &
\multicolumn{1}{|l|}{$0.122929$} &
\multicolumn{1}{|l|}{$0.148044$} &
\multicolumn{1}{|l|}{$0.098305$} &
\multicolumn{1}{|l|}{$0.080768$} \\ \hline
\multicolumn{1}{|l|}{$34a1$} & \multicolumn{1}{|l|}{$83$} &
\multicolumn{1}{|l|}{$0.130011$} &
\multicolumn{1}{|l|}{$0.103015$} &
\multicolumn{1}{|l|}{$0.141567$} &
\multicolumn{1}{|l|}{$0.085853$} &
\multicolumn{1}{|l|}{$0.077536$} \\ \hline
\multicolumn{1}{|l|}{$34a1$} & \multicolumn{1}{|l|}{$89$} &
\multicolumn{1}{|l|}{$0.132737$} &
\multicolumn{1}{|l|}{$0.115162$} &
\multicolumn{1}{|l|}{$0.150665$} & \multicolumn{1}{|l|}{$0.0967$}
& \multicolumn{1}{|l|}{$0.08984$} \\ \hline
\multicolumn{1}{|l|}{$34a1$} & \multicolumn{1}{|l|}{$93$} &
\multicolumn{1}{|l|}{$0.127371$} &
\multicolumn{1}{|l|}{$0.120887$} &
\multicolumn{1}{|l|}{$0.147903$} &
\multicolumn{1}{|l|}{$0.100794$} &
\multicolumn{1}{|l|}{$0.071056$} \\ \hline
\multicolumn{1}{|l|}{$34a1$} & \multicolumn{1}{|l|}{$101$} &
\multicolumn{1}{|l|}{$0.124347$} &
\multicolumn{1}{|l|}{$0.118514$} &
\multicolumn{1}{|l|}{$0.136806$} & \multicolumn{1}{|l|}{$0.10392$}
& \multicolumn{1}{|l|}{$0.079833$} \\ \hline
\multicolumn{1}{|l|}{$34a1$} & \multicolumn{1}{|l|}{$115$} &
\multicolumn{1}{|l|}{$0.13509$} & \multicolumn{1}{|l|}{$0.116645$}
& \multicolumn{1}{|l|}{$0.140947$} &
\multicolumn{1}{|l|}{$0.089575$} &
\multicolumn{1}{|l|}{$0.075199$} \\ \hline
\multicolumn{1}{|l|}{$34a1$} & \multicolumn{1}{|l|}{$117$} &
\multicolumn{1}{|l|}{$0.130015$} &
\multicolumn{1}{|l|}{$0.124739$} &
\multicolumn{1}{|l|}{$0.140575$} &
\multicolumn{1}{|l|}{$0.104611$} &
\multicolumn{1}{|l|}{$0.078912$} \\ \hline
\multicolumn{1}{|l|}{$34a1$} & \multicolumn{1}{|l|}{$123$} &
\multicolumn{1}{|l|}{$0.135871$} &
\multicolumn{1}{|l|}{$0.107328$} &
\multicolumn{1}{|l|}{$0.133703$} &
\multicolumn{1}{|l|}{$0.093233$} & \multicolumn{1}{|l|}{$0.07586$}
\\ \hline
\end{tabular}

\begin{tabular}{|c|c|c|c|c|c|c|}
\hline $E$ & $n_{0}$ & $225$ & $256$ & $289$ & $324$ & $361$ \\
\hline \multicolumn{1}{|l|}{$34a1$} & \multicolumn{1}{|l|}{$1$} &
\multicolumn{1}{|l|}{$0.120486$} & \multicolumn{1}{|l|}{$0.06472$}
& \multicolumn{1}{|l|}{$0.062951$} &
\multicolumn{1}{|l|}{$0.070443$} &
\multicolumn{1}{|l|}{$0.055249$} \\ \hline
\multicolumn{1}{|l|}{$34a1$} & \multicolumn{1}{|l|}{$13$} &
\multicolumn{1}{|l|}{$0.120107$} & \multicolumn{1}{|l|}{$0.0621$}
& \multicolumn{1}{|l|}{$0.0692$} &
\multicolumn{1}{|l|}{$0.077468$} &
\multicolumn{1}{|l|}{$0.055683$} \\ \hline
\multicolumn{1}{|l|}{$34a1$} & \multicolumn{1}{|l|}{$19$} &
\multicolumn{1}{|l|}{$0.077919$} &
\multicolumn{1}{|l|}{$0.059063$} &
\multicolumn{1}{|l|}{$0.053673$} &
\multicolumn{1}{|l|}{$0.070172$} &
\multicolumn{1}{|l|}{$0.038649$} \\ \hline
\multicolumn{1}{|l|}{$34a1$} & \multicolumn{1}{|l|}{$21$} &
\multicolumn{1}{|l|}{$0.118898$} &
\multicolumn{1}{|l|}{$0.062197$} & \multicolumn{1}{|l|}{$0.06578$}
& \multicolumn{1}{|l|}{$0.071127$} &
\multicolumn{1}{|l|}{$0.056771$} \\ \hline
\multicolumn{1}{|l|}{$34a1$} & \multicolumn{1}{|l|}{$33$} &
\multicolumn{1}{|l|}{$0.118185$} &
\multicolumn{1}{|l|}{$0.059254$} &
\multicolumn{1}{|l|}{$0.065998$} &
\multicolumn{1}{|l|}{$0.072214$} &
\multicolumn{1}{|l|}{$0.053245$} \\ \hline
\multicolumn{1}{|l|}{$34a1$} & \multicolumn{1}{|l|}{$35$} &
\multicolumn{1}{|l|}{$0.092942$} &
\multicolumn{1}{|l|}{$0.058685$} &
\multicolumn{1}{|l|}{$0.050621$} &
\multicolumn{1}{|l|}{$0.063069$} &
\multicolumn{1}{|l|}{$0.039964$} \\ \hline
\multicolumn{1}{|l|}{$34a1$} & \multicolumn{1}{|l|}{$43$} &
\multicolumn{1}{|l|}{$0.097305$} &
\multicolumn{1}{|l|}{$0.057743$} &
\multicolumn{1}{|l|}{$0.053748$} &
\multicolumn{1}{|l|}{$0.069417$} & \multicolumn{1}{|l|}{$0.0381$}
\\ \hline \multicolumn{1}{|l|}{$34a1$} &
\multicolumn{1}{|l|}{$53$} & \multicolumn{1}{|l|}{$0.120493$} &
\multicolumn{1}{|l|}{$0.060345$} &
\multicolumn{1}{|l|}{$0.065208$} &
\multicolumn{1}{|l|}{$0.071921$} &
\multicolumn{1}{|l|}{$0.049442$} \\ \hline
\multicolumn{1}{|l|}{$34a1$} & \multicolumn{1}{|l|}{$59$} &
\multicolumn{1}{|l|}{$0.097042$} &
\multicolumn{1}{|l|}{$0.061403$} &
\multicolumn{1}{|l|}{$0.047947$} &
\multicolumn{1}{|l|}{$0.065448$} &
\multicolumn{1}{|l|}{$0.040482$} \\ \hline
\multicolumn{1}{|l|}{$34a1$} & \multicolumn{1}{|l|}{$67$} &
\multicolumn{1}{|l|}{$0.09381$} & \multicolumn{1}{|l|}{$0.057656$}
& \multicolumn{1}{|l|}{$0.066198$} &
\multicolumn{1}{|l|}{$0.064485$} &
\multicolumn{1}{|l|}{$0.038135$} \\ \hline
\multicolumn{1}{|l|}{$34a1$} & \multicolumn{1}{|l|}{$69$} &
\multicolumn{1}{|l|}{$0.121028$} &
\multicolumn{1}{|l|}{$0.065232$} &
\multicolumn{1}{|l|}{$0.066193$} &
\multicolumn{1}{|l|}{$0.066645$} &
\multicolumn{1}{|l|}{$0.052697$} \\ \hline
\multicolumn{1}{|l|}{$34a1$} & \multicolumn{1}{|l|}{$77$} &
\multicolumn{1}{|l|}{$0.120666$} &
\multicolumn{1}{|l|}{$0.065448$} &
\multicolumn{1}{|l|}{$0.067139$} &
\multicolumn{1}{|l|}{$0.072256$} &
\multicolumn{1}{|l|}{$0.052677$} \\ \hline
\multicolumn{1}{|l|}{$34a1$} & \multicolumn{1}{|l|}{$83$} &
\multicolumn{1}{|l|}{$0.098931$} &
\multicolumn{1}{|l|}{$0.058119$} &
\multicolumn{1}{|l|}{$0.049453$} &
\multicolumn{1}{|l|}{$0.067301$} &
\multicolumn{1}{|l|}{$0.040662$} \\ \hline
\multicolumn{1}{|l|}{$34a1$} & \multicolumn{1}{|l|}{$89$} &
\multicolumn{1}{|l|}{$0.117926$} &
\multicolumn{1}{|l|}{$0.059319$} &
\multicolumn{1}{|l|}{$0.062727$} &
\multicolumn{1}{|l|}{$0.072717$} &
\multicolumn{1}{|l|}{$0.055294$} \\ \hline
\multicolumn{1}{|l|}{$34a1$} & \multicolumn{1}{|l|}{$93$} &
\multicolumn{1}{|l|}{$0.115051$} &
\multicolumn{1}{|l|}{$0.060207$} &
\multicolumn{1}{|l|}{$0.061973$} &
\multicolumn{1}{|l|}{$0.076746$} &
\multicolumn{1}{|l|}{$0.058861$} \\ \hline
\multicolumn{1}{|l|}{$34a1$} & \multicolumn{1}{|l|}{$101$} &
\multicolumn{1}{|l|}{$0.123836$} &
\multicolumn{1}{|l|}{$0.063449$} &
\multicolumn{1}{|l|}{$0.070854$} &
\multicolumn{1}{|l|}{$0.072985$} &
\multicolumn{1}{|l|}{$0.058416$} \\ \hline
\multicolumn{1}{|l|}{$34a1$} & \multicolumn{1}{|l|}{$115$} &
\multicolumn{1}{|l|}{$0.099212$} &
\multicolumn{1}{|l|}{$0.057176$} &
\multicolumn{1}{|l|}{$0.049053$} &
\multicolumn{1}{|l|}{$0.065697$} &
\multicolumn{1}{|l|}{$0.037937$} \\ \hline
\multicolumn{1}{|l|}{$34a1$} & \multicolumn{1}{|l|}{$117$} &
\multicolumn{1}{|l|}{$0.120381$} &
\multicolumn{1}{|l|}{$0.061819$} &
\multicolumn{1}{|l|}{$0.067332$} & \multicolumn{1}{|l|}{$0.07149$}
& \multicolumn{1}{|l|}{$0.054083$} \\ \hline
\multicolumn{1}{|l|}{$34a1$} & \multicolumn{1}{|l|}{$123$} &
\multicolumn{1}{|l|}{$0.097763$} &
\multicolumn{1}{|l|}{$0.057011$} &
\multicolumn{1}{|l|}{$0.049838$} & \multicolumn{1}{|l|}{$0.06636$}
& \multicolumn{1}{|l|}{$0.038657$} \\ \hline
\end{tabular}
\end{center}

\subsection{A Graphical Example\protect\bigskip}

We plot some graph of the data for $E=11a1$, $n_{0}=1$, $k=1$. In
this graph on the $x-$axis we plot $x_{n_{0}}(M)$ up to
$M=10^{7}$. Dots above at the beginning belong to the graph of the
function $s_{1,1}(x_{1}(M))/x_{1}(M)$,
dots below at the beginning belong to the graph of the function $0.295669%
\frac{(\log \log (x_{1}(M)))^{1.005}}{\log (x_{1}(M))}.$

Received: 27 September, 2010 and in revised form 24 January 2011.
\end{document}